\documentclass [a4paper,notitlepage,11pt]{article}
\usepackage[pdftex]{graphicx}
\usepackage{graphicx}
\usepackage{float}
\usepackage{color}
\usepackage{hyperref}
\hypersetup{colorlinks=true,
            linkcolor={blue},
            citecolor={green}}
\usepackage[shortlabels]{enumitem}

\usepackage{amsmath, amsfonts, latexsym}
\usepackage{mathabx}
\usepackage{mathrsfs,amssymb}
\usepackage{version}
\usepackage{bm}
\usepackage{bbm}
\usepackage{flafter}

\usepackage{dsfont}
\usepackage{amsthm}
\usepackage[english]{babel}

\numberwithin{equation}{section}

\newcommand{\E}{\mathbb E}

\newcommand{\R}{\mathbb R}
\newcommand{\p}{\mathbb P}
\newcommand{\Q}{\mathbb Q}

\newcommand{\M}{\mathbb M}

\newcommand{\1}{\mathbbm 1}
\newcommand{\defeq}{\mathrel{\mathop:}=}

\usepackage{bbm}

\DeclareMathOperator{\Tr}{Tr}

\newcommand{\ti}[1]{\ensuremath{ t_{ #1} }}
\newcommand{\esp}[1]
    {\ensuremath{%
     \mathbb{E}\!\!\left[#1\right]}}

 \theoremstyle{plain}
 \newtheorem{theo}{Theorem}[section]
 \newtheorem{lem}{Lemma}[section]
 \newtheorem{prop}{Proposition}[section]

 \theoremstyle{definition}
 \newtheorem{dfn}{Definition}[section]

 \def \proof{{\noindent \bf Proof. }}
 \def \eproof{\hbox{ }\hfill$\Box$}

 \theoremstyle{remark}
 \newtheorem{rem}{Remark}[section]

  \newtheorem*{hypot}{Hypothesis}

\newcommand{\set}[1]
    {\ensuremath{\left \{ #1 \right \}}}

\renewcommand{\ti}[1]
    {\ensuremath{t_{i#1}}}

\newcommand{\bT}
    {\ensuremath{ \mathbb{T}}}

  \def\P{\mathbb{P}}
 \def\Fc{\mathcal{F}}
  \def\Ec{\mathcal{E}}
    \def\Uc{\mathcal{U}}
      \def\Ac{\mathcal{A}}
  \def\d{\mathrm{d}}
  \def\1{\mathbf{1}}
  \def\S{\mathrm{S}}
  \def\veee{\!\vee\!}
  \def\x{\times}
  \def\ti{t_{i}}
    \def\tip{t_{i+1}}
      
      \def\eps{\varepsilon}
       
       \def\Tr#1{{\rm Tr}\left[#1\right]}
        \def\Hi{{\bf${\rm H_{i}}$}}
 \def\Hip{{\bf${\rm H_{i+1}}$}}
 \def\co#1{{\rm co}[#1]}

 \def\Rpp{{\mathcal O}_{\!+}^{d}}
 \def\bru#1{{\color{red}#1}}

  \renewcommand{\ell}{g}
   \renewcommand{\zeta}{\tilde{\ell}}
  
 \title{
A backward dual representation for the quantile hedging of Bermudan options} 
 \author{Bruno Bouchard\thanks{Research supported by ANR Liquirisk and Investissements d'Avenir (ANR-11-IDEX-0003/Labex Ecodec/ANR-11-LABX-0047).}
         \\\small CEREMADE, Universit\'{e} Paris Dauphine
         \\\small and CREST-ENSAE
         \\\small \sf bouchard@ceremade.dauphine.fr
         \and
         G\'{e}raldine Bouveret\thanks{Research supported by the Natixis Foundation for Quantitative Research.}
         \\\small Department of Mathematics
         \\\small Imperial College London
         \\\small \sf g.bouveret11@imperial.ac.uk
         \and
         {Jean-Fran\c{c}ois Chassagneux}\footnotemark[1]
         \\{\small Laboratoire  de Probabilit\'{e}s  et Mod\`{e}les  Al\'{e}atoires}
         \\{\small CNRS, UMR 7599,  Universit\'{e} Paris Diderot}
         \\{\small \sf jean-francois.chassagneux@univ-paris-diderot.fr}
}

\usepackage{ulem}

\begin{document}

 \maketitle

 \begin{abstract} \noindent Within a Markovian complete financial market, we consider the problem of hedging a Bermudan option with a given probability. Using stochastic target and duality arguments, we derive a backward algorithm for the Fenchel transform of the pricing function. This algorithm is similar to the usual American backward induction, except that it requires two additional Fenchel transformations at each exercise date. We provide numerical illustrations.
    \end{abstract}


 \vspace{1em}


\section{Introduction}

We study the problem of hedging a claim of Bermudan style with a given probability $p$. More precisely, we want to characterize the minimal initial value $v(\cdot,p)$ of an hedging portfolio for which we can find a financial strategy such that, with a probability $p$, it  remains above the exercise value of the Bermudan option at any possible exercise date.

This problem is referred to {\sl quantile hedging},  and it was popularized by   F\"ollmer et al.~\cite{follmer1999quantile,follmer2000efficient} . For  claims of European type, they explained how   the so-called {\sl quantile hedging price} can be computed explicitly when the market is complete,    by using duality arguments or the Neyman-Pearson lemma. A similar question was studied in Bouchard et al.~\cite{bouchard2009stochastic} but in a Markovian setting. They showed that, even in incomplete markets and for general {\sl loss functions}, one can characterize the pricing function as the solution of a non-linear parabolic second order differential equation, by using tools developed in the context of stochastic target problems by Soner and Touzi \cite{soner2002stochastic}-\cite{soner2003stochastic}. When the market is complete, they also observed that taking a Legendre-Fenchel transform in the equation reduces the computation of the price to the resolution of a linear parabolic second order differential equation, which can be solved {\sl explicitly} by using the Feynman-Kac formula.

As far as super-hedging is concerned, the pricing of a  Bermudan option reduces to a backward sequence of pricing problems for European claims. It is therefore natural to ask whether a similar result holds for the quantile hedging price, and whether one can   extend the   closed-form solutions of \cite{follmer1999quantile} and   \cite{bouchard2009stochastic}  to Bermudan options.

This  paper answers to the positive.   Namely, we provide a backward induction algorithm for the Fenchel transform $w$ of the quantile hedging price $v(\cdot,p)$, with respect to the parameter $p$ which prescribes the probability of hedging, see \eqref{eq: def w} and Theorem \ref{th mainrsult}. The algorithm \eqref{eq: def w}  is in a sense very similar to the one used for the pricing of Bermudan options. It is however written on the  Fenchel transform $w$, rather than $v$, and  it involves two additional Fenchel transformations at each exercise date.

To derive this, we  first build on the original idea of \cite{bouchard2009stochastic} which consists
in increasing the state space in order to reduce to a stochastic target problem of American type,
as studied in Bouchard and Vu \cite{bouchard2010obstacle}. We then follow a very different route.
Instead of appealing to stochastic target technics, we  derive from this formulation
a first dynamic programming algorithm for  $v$, see Proposition \ref{wProbRep0}, which
relates to a series of  optimal control of martingale problems.
This dynamic programming principle suggests  a backward algorithm for the computation of the Fenchel transform. It is  defined in  \eqref{eq: def w}. We  analyze it in details  in Section \ref{sec: travaille sur w}. The main difficulty consists in controlling the propagation of the differentiability and growth properties  of the corresponding value function, backward in time.  Then, as in \cite{bouchard2009stochastic,follmer1999quantile}, a martingale representation argument allows us to show, by backward induction, that the algorithm in \eqref{eq: def w} and Proposition \ref{wProbRep0} provides the Fenchel transform of one another.

Before concluding this introduction, we would like to point out  that a similar
problem has been studied recently by   Jiao et al.~\cite{jiao2013hedging} in the
 form of general {\sl lookback-style contraints}. They provide an alternative
formulation in terms of an optimal control of martingale problems. This has to
be compared with \cite{reveillac2012bsdes} and our Proposition \ref{wProbRep0}.
No Markovian structure is required, but  they do not provide an explicit scheme
as we do. Moreover, the smoothness conditions they impose on their loss
functions are not satisfied in the quantile hedging case.
They also study the case of several constraints in expectation set (independently) at the different exercise times, which is close to the P\&L matching problems of   Bouchard and Vu \cite{bouchard2012stochastic}.

Finally, in this paper, we focus on the quantile hedging problem for sake of simplicity.
It is an archetype of an irregular loss function, and it should be clear
that a similar analysis can be carried out for a wide class of  (more regular) loss functions.
Also note that to obtain the dual algorithm, we only use probabilistic arguments 
which opens the door to the study of more general non-Markovian settings. 
\vspace{2mm}

   \textbf{Notations:}
   Let $d$ be a positive integer. Any vector $x$ of $\R^{d}$ is seen as a
column vector. Its norm and transpose are denoted by $|x|$ and $x^{\top}$. We
set $\mathbb{M}^{d}:=\mathbb{R}^{d\x d}$ and denote by $M^{\top}$ the transpose
of $M\in \mathbb{M}^{d}$, while $\Tr M$ is its trace. For ease of notations, we
set $\Rpp:=(0,\infty)^{d}$.

 We fix a finite time horizon $T>0$.  Let $\psi : (t,x,p)\in [0,T]\x\Rpp\x \R
\mapsto \psi(t,x,p)$. If it is smooth enough, we denote by $\partial_{t}\psi$
and $\partial_{p}\psi$   its   derivative with respect $t$ and $p$, and by
$\partial_{x} \psi$ its Jacobian matrix with respect to $x$, as a column vector.
The Hessian with respect to $x$ is $\partial^{2}_{xx} \psi$,
$\partial^{2}_{pp}\psi$ is the second order derivative with respect to $p$, and
$\partial^{2}_{xp} \psi$ is the vector of cross second order derivatives. We
denote by $\psi^{\sharp}$ its Fenchel transform with respect to the last
argument,
\begin{eqnarray}\label{eq :def Fenchel transform}
 \psi^{\sharp}(t,x,q):=\sup_{p\in \R}\left(pq - \psi(t,x,p)\right)\,,
 \end{eqnarray}
 and define
 \begin{eqnarray*}
\co\psi \;,\text{the closed convex envelope of $\psi$ with respect to its
last argument.}
 \end{eqnarray*}
 If $\psi$ is convex with respect to its last variable, we denote by $\mathrm D^{+}_{p}\psi$ and $\mathrm D^{-}_{p}\psi$ its corresponding right- and left-derivatives. We refer to \cite{rocka} for the various notions related to convex analysis.

  We fix a complete probability space $(\Omega,\mathbb{F},\p)$  supporting a $d$-dimensional Brownian motion
$W$. We denote by $\mathbb{F} = (\mathcal{F}_t)_{0 \le t \le T}$ the usual augmented  Brownian
filtration.  All over the paper,  inequalities between random variables have to be understood in the $\p$-a.s. sense.

\setlength{\parindent}{0cm}
\section{Problem formulation and main results}

\subsection{Financial market and hedging problem}

  Our  financial market consists in a non-risky
asset, whose price process is normalized to unity, and $d$
risky assets $X=(X^1,...,X^d)$ whose dynamics are given by
    \begin{align}\label{X}
      X^{t,x}_{s}= x + \int_t^s \mu(r, X^{t,x}_{r}) \mathrm{d} r +
\int_t^s \sigma(r,X^{t,x}_{r})\mathrm{d} W_r\,,
    \end{align}
given the initial data $(t,x)\in[0,T]\times\Rpp$. To ensure that the above is   well-defined, we assume that
\begin{eqnarray}\label{eq: cond mu sig}
\mu:[0,T]\x{\Rpp}\rightarrow\R^d \mbox{ and
}\sigma:[0,T]\x{\Rpp}\rightarrow\M^d \mbox{ are Lipschitz continuous\,,}
\end{eqnarray}
 and that
the unique strong  solution to \eqref{X} takes its values in $\Rpp$ when the original data lies in $\Rpp$.

In order to enforce the absence of arbitrage  and the completeness of the financial market, we also impose that
\begin{align}\label{eq: cond lambda}
       &\sigma \text{ is invertible\,, }   \;\;
         \lambda\defeq\sigma^{-1}\mu \,\text{ is bounded}\\
         &\text{ and  Lipschitz continuous in space, uniformly in time}\,.\nonumber
    \end{align}
The Lipschitz continuity condition is not required to define the risk neutral measure\footnote{$\Ec$ denotes here the Dol\'eans-Dade exponential.}
\begin{eqnarray}\label{eq: def Q}
\Q_{t,x}:=\frac{1}{Q^{t,x,1}_{T}}\cdot \P\; \mbox{ with }\; \frac{1}{Q^{t,x,q}}:=\frac{1}{q}\Ec\left(-\int_{t}^{\cdot} \lambda(s,X^{t,x}_{s})^{\top} \d W_{s} \right)\,,\;q>0\,,
\end{eqnarray}
but will be used in some of our forthcoming arguments.

In this model, an admissible financial strategy is a  $d$-dimensional   predictable process $\nu$ such that
    \begin{equation}\label{eq: cond carre int}
      \E^{\Q_{t,x}}\left[\int_{t}^{T} \|\nu_{s}^{\top}\sigma(s,X^{t,x}_{s})\|^{2}ds\right]
      <\infty\,,
    \end{equation}
and the corresponding wealth process remains non-negative
    \begin{equation*}
      Y^{t,x,y,\nu}:=y+\int_{t}^{\cdot} \nu_{r}^{\top}  \mathrm{d}X^{t,x}_{r} \ge 0\,, \;\text{ on } [t,T]\,,
    \end{equation*}
   given the initial data $(t,x)$ of the market and the initial dotation $y\ge 0$. We denote by $\Uc_{t,x,y}$ the collection of admissible financial strategies.
As usual, each $\nu^{i}_{t}$ should be interpreted as the number of units of asset $i$ in the portfolio at time $t$.

We now fix a finite collection of times
    \begin{align*}
      \bT_{t}\defeq \set{t_0=0\leq \dots \leq t_i\leq \dots\leq t_n=T}\cap (t,T]\,,
    \end{align*}
together with payoff functions
\begin{eqnarray}\label{eq: ass ell lipschitz}
x\in \Rpp \mapsto \ell(t_i,x)\ge 0,\;\mbox{Lipschitz continuous for all
$i\le n$}\,.
\end{eqnarray}

Our quantile hedging problem consists in finding the minimal initial wealth
$v(t,x,p)$ which ensures that the stream of Bermudan payoffs $\{
\ell(s,X^{t,x}_{s}),\;s\in {\bT_{t}}\}$ can be hedged with a given
probability  $p$, \begin{eqnarray}\label{eq: def v}
v(t,x,p):=\inf\Gamma(t,x,p)\,,
\end{eqnarray}
where
\begin{eqnarray*}
&\Gamma(t,x,p):=\left\{y\ge 0~:~\exists\, \nu \in \Uc_{t,x,y}\mbox{ s.t. }\p\left[\bigcap_{s\in \bT_{t}} \S^{t,x,y,\nu}_{s}\right]\ge p\right\}\;,& \\
&\mbox{ with } \;  {
 \S^{t,x,y,\nu}_{s}:=\left\{
 \begin{array}{lcl}
 \Omega & \mbox{if }&    s \le t \\
 \{Y^{t,x,y,\nu}_{s}\ge \ell(s,X^{t,x}_{s})\}   & \mbox{if }& s>t
\end{array}\right.}.  &
\end{eqnarray*}

Observe that $v(t,\cdot)$ must be  interpreted as a
{continuation}  value, i.e.~the price at time $t$ knowing that the option
has not been exercised on $[0,t]$. In particular,   $v(T,\cdot)=0$. For $p=1$,
$v(t,\cdot,1)$ coincides with the continuation value of the super-hedging price
of the Bermudan option. In this complete market, it satisfies the usual dynamic
programming principle
\begin{eqnarray}\label{eq: dpp retro pour v(t,x,1)}
v(t,x,1)&=& \E^{\Q_{t,x}}[(v\veee \ell)(t_{i+1},X^{t,x}_{t_{i+1}},1)]\,,\; \text{ for } t\in [t_{i},t_{i+1})\,,\;i<n\,,
\end{eqnarray}
{see \cite{schweizer2002bermudan}.}
Above and in the following, we use the notation
$$
 \ell(t,x,p):=\ell(t,x)\1_{\{0<p\le 1\}} +\infty \1_{\{p>1\}}\,, \; { \text{
for } p \in \R}\,.
$$

Note that $\Gamma$ can also be formulated in terms of stopping times, see the Appendix for the proof.
\begin{prop}\label{prop : formulation sur temps d'adrets}   For $(t,x,p)\in [0,T]\x \Rpp\x [0,1]$,
\begin{eqnarray}\label{eq: formulation Gamma avec temps arret}
\begin{array}{rcl}\Gamma(t,x,p)&=& \{y\ge 0~:~\exists\, \nu \in
\Uc_{t,x,y}{\;\;\mbox{\rm s.t. }}\p [ \S^{t,x,y,\nu}_{\tau} ]\ge
p,\;\forall\;\tau \in {\cal T}_{t} \}\\
&=&
\{y\ge 0~:~\exists\, \nu \in \Uc_{t,x,y}{\;\;\mbox{\rm s.t. }}\p [
\S^{t,x,y,\nu}_{\hat \tau_{\nu}} ]\ge p \}\1_{\{t<T\}}+\R_{+}\1_{\{t=T\}}\,,
\end{array}
\end{eqnarray}
in which ${\cal T}_{t}$ is the set of stopping times with values in $\bT_{t}$, and $\hat  \tau_{\nu}:=\min\{s\in \bT_{t}~:~Y^{t,x,y,\nu}_{s}<\ell(s,X^{t,x}_{s})\}\wedge T$.
 \end{prop}

\begin{rem}\label{rem : monotony and pmin}  The function $p\mapsto v(\cdot,p)$ is   non-decreasing.
It takes the value $0$ if $p\le p_{\rm min}(t,x)$ where
\begin{eqnarray}\label{eq: def pmin}
p_{\rm min}(t,x):=\P[\ell(s,X^{t,x}_{s})=0\;\text{ for all } s\in \bT_{t}]\,,
\end{eqnarray}
with the convention $p_{\rm min}(T,\cdot)=1$.
To avoid trivial statements, we assume that $p_{\rm min}(t,\cdot)<1$, for $t<T$, which implies
\begin{align}
v(t,x,1) > 0\,,\; \text{for} \;t < T\,. \label{eq v(t,x,1) strict pos}
\end{align}

Moreover, it follows from \eqref{eq: ass ell lipschitz}   that we can find $C>0$
such that $\ell(s,x)\le C(1+\sum_{i=1}^{d} x^{i})$, for $x\in \Rpp$, $s \in
\bT_0$. This implies that we can restrict to strategies $\nu$ such that
\begin{eqnarray}\label{eq: borne Y}
0\le Y^{t,x,y,\nu}\le C(1+ |X^{t,x}|)\,,
\end{eqnarray}
{by possibly adopting a buy-and-hold strategy after the first time at which the
wealth process hits the right-hand side term, recall that $X^{t,x}$ has positive
components.} In particular,
\begin{eqnarray}\label{eq: borne v}
0\le v(t,x,p)\le C(1+ |x|)\,.
\end{eqnarray}
\end{rem}

\subsection{Equivalent formulation as a stochastic target problem}

The first step in our analysis consists in  reducing the problem to a stochastic target problem of American type as studied in  \cite{bouchard2010obstacle}. As in \cite{bouchard2009stochastic}, we first increase the dimension of the controlled process by introducing the family of martingales
$$
P^{t,p,\alpha}:=p+\int_{t}^{\cdot}\alpha_{s}^{\top}\d W_{s}\,,
$$
where $\alpha$ is a square integrable predictable process. The process
$P^{t,p,\alpha}$ will be later on interpreted as the conditional probability of
success. It is therefore natural to restrict to the class of controls  such that
$$
P^{t,p,\alpha}\in [0,1]\,,\;\text{ on }\; [t,T]\,.
$$

We denote by $\Ac_{t,p}$ the set of predictable square integrable processes such
that the above holds, and set $\hat \Uc_{t,x,y,p}:=\Uc_{t,x,y}\times \Ac_{t,p}$.

\begin{prop} \label{pr PbReduction0}      Fix   $(t,x,p) \in [0,T]\x \Rpp\times[0,1]$, then
      \begin{gather}\label{eq: Gamma = formulation augmentee}
        {\Gamma}(t,x,p)
       = \Big\{y\geq 0: \exists \;(\nu,\alpha)\in  \hat \Uc_{t,x,y,p}\; \text{ s.t. }\;Y^{t,x,y,\nu} \geq \ell(\cdot,X^{t,x},P^{t,p,\alpha})\;\text{ on } \bT_{t}\Big\}\,.
      \end{gather}
    \end{prop}
 \proof At time $T$ both sets are $\R_{+}$ by definition of $\bT_{T}$. We now fix $t<T$.   Let $\bar \Gamma(t,x,p)$ denote the right-hand side in \eqref{eq: Gamma = formulation augmentee} and let $y$ be one of his elements.
  Fix $(\nu,\alpha)\in \hat \Uc_{t,x,y,p}$ such that $Y^{t,x,y,\nu} \geq \ell(\cdot,X^{t,x},P^{t,p,\alpha})$ on $\bT_{t}$. Then, $\S^{t,x,y,\nu} \supset \{P^{t,p,\alpha} >0\}$ on $\bT_{t}$. Since   $P^{t,p,\alpha}\in [0,1]$ and therefore  $\1_{\{P^{t,p,\alpha} >0\} }\ge P^{t,p,\alpha}$, this implies
 \begin{eqnarray*}
 \p\left[\cap_{s\in \bT_{t}} \S^{t,x,y,\nu}_{s} \right]&\ge& \p\left[\cap_{s\in \bT_{t}} \{P^{t,p,\alpha}_{s} >0\} \right]
   \ge  \E\left[P^{t,p,\alpha}_{T}\prod_{s\in \bT_{t}\setminus\{T\}}  \1_{\{P^{t,p,\alpha}_{s} >0\} }\right]\,.
 \end{eqnarray*}
  The process $P^{t,p,\alpha}$ being a martingale, $\{P^{t,p,\alpha}_{s}=0\}\subset  \{P^{t,p,\alpha}_{T}=0\}$,  {$s\in(t,T]$}. Hence
  \begin{eqnarray*}
 \p\left[\cap_{s\in \bT_{t}} \S^{t,x,y,\nu}_{s} \right]&\ge& \E\left[P^{t,p,\alpha}_{T} \right]=p\,.
 \end{eqnarray*}
 Therefore, $y\in \Gamma(t,x,p)$ and this argument proves that
 $\bar \Gamma(t,x,p)\subset  \Gamma(t,x,p)$.

 We now fix $y\in   \Gamma(t,x,p)$ and choose $\nu \in \Uc_{t,x,y}$ such that $p':=\p\left[\bigcap_{s\in \bT_{t}} \S^{t,x,y,\nu}_{s}\right]\ge p$. By the martingale representation theorem, we can find $\alpha\in \Ac_{t,p'}$ such that
 $$
 \1_{\bigcap_{s\in \bT_{t}} \S^{t,x,y,\nu}_{s}}=P^{t,p',\alpha}_{T}\ge P^{t,p,\alpha}_{T}\,.
 $$
 By possibly replacing $\alpha$ by the constant process $0$ after the first time after $t$ at which  $P^{t,p,\alpha}$ reaches the level $0$, we can assume that $\alpha\in \Ac_{t,p}$. Moreover, the above implies
 $$
 \1_{ \S^{t,x,y,\nu}_{s}} \ge P^{t,p,\alpha}_{T}\,,\;s\in \bT_{t}\,,
$$
which by taking the conditional expectation and using the fact that $ P^{t,p,\alpha}$ is a martingale leads to
 $
 \1_{ \S^{t,x,y,\nu} }\ge P^{t,p,\alpha}$   on $\bT_{t}$. The latter is equivalent to $Y^{t,x,y,\nu} \geq \ell(\cdot,X^{t,x},P^{t,p,\alpha})$   on $\bT_{t}$.
Hence, $y\in\bar \Gamma(t,x,p)$.
\eproof

\subsection{Dynamic programming and dual backward algorithm}\label{se DPP}

With the formulation obtained in Proposition \ref{pr PbReduction0} at hand, one can now derive a first dynamic programming algorithm. Its proof is postponed to the Appendix.
 \begin{prop} \label{wProbRep0}    Fix $0\leq i\leq n-1$  and
$(t,x,p)\in[t_i,t_{i+1})\x\Rpp\times[0,1]$,
\begin{align}
v(t,x,p)&=\inf_{\alpha\in\mathcal{A}_{t,p}}\E^{\Q_{t,x}}\left[(v\veee \ell)\left(t_{i+1},X^{ t ,x}_{t_{i+1}},P^{t,p,\alpha}_{t_{i+1}}\right)\right]\,.
        \label{eq probw0}
\end{align}
 As a consequence, there exists $C>0$ such that
\begin{eqnarray}\label{eq: v loc lip}
|v(t,x,p)-v(t,x',p)|\le C(1+|x|+|x'|)|x-x'|\,,
\end{eqnarray}
for all  $(t,p)\in [0,T]\x [0,1]$ and $x,x'\in \Rpp$.
\end{prop}

\begin{rem}\label{rem convexifcation primal} We shall see in Section \ref{sec: proofs} that $(v\veee \ell)$ can be replaced by its convex envelope with respect to $p$ in  \eqref{eq probw0}. This phenomenon was already observed in \cite{reveillac2012bsdes} and \cite{bouchard2009stochastic}.
\end{rem}

Note that this provides a first way to compute the value function $v$. Indeed,
standard arguments (see \cite{bouchard2009stochastic}) should lead to a characterization of  $v$ on each interval $[\ti,\tip)$, $i<n$ and on $\Rpp\times(0,1)$ as a viscosity
solution of
\begin{align}\label{eq: Hjb1}
\sup_{\alpha\in \R^{d}} \left\{-\partial_t\varphi(\cdot) + \alpha^\top\lambda \partial_p\varphi(\cdot) -\frac12\left(  \Tr{\sigma\sigma^{\top} \partial^2_{xx}\varphi(\cdot)}+2\alpha^{\top}\sigma^{\top} \partial^2_{xp}\varphi(\cdot)
+|\alpha|^{2}\partial^2_{pp}\varphi(\cdot)\right) \right\}{=0}\,,
\end{align}
with the boundary conditions
\begin{gather}
v(\tip-,\cdot)=(v\veee \ell)(\tip,\cdot),\,\text{on}\,\Rpp\times[0,1]\label{eq: Hjb bord1}\\
v(\cdot,1)=\E^{\Q_{t,x}}[(v\veee \ell)(t_{i+1},X^{t,x}_{t_{i+1}},1)],\quad v(\cdot,0)=0,\,\text{on}\,[\ti,\tip)\times\Rpp,\,i<n\,.\label{eq: Hjb bord2}
\end{gather}
However, the fact that the control $\alpha\in \R^{d}$ in the above is not bounded (as it comes from the martingale representation theorem) makes the associated Hamilton-Jacobi-Bellman operator in \eqref{eq: Hjb1} discontinuous. More precisely it is lower semi-continuous but not upper semi-continuous and a precise statement would then require a relaxation of the operator in \eqref{eq: Hjb1}. 
This discontinuity makes the proof of a comparison result very difficult and the latter is necessary to build convergent numerical schemes. One way to overcome this problem is to consider instead the Fenchel transform $v^{\sharp}$ of $v$, see \eqref{eq :def Fenchel transform} in the notations section.

 Indeed, heuristically, as already observed in \cite{bouchard2009stochastic} in the case $n=1$, a change of variable argument in \eqref{eq: Hjb1} and the exploitation of the boundary conditions in \eqref{eq: Hjb bord2} suggests that  the dual function $v^{\sharp}$ should be at least {a viscosity} sub-solution of the linear partial differential equation
  \begin{gather}
         -\partial_t\varphi(\cdot) -\frac{1}{2}\left({\rm Tr}[\sigma\sigma^{\top} \partial^{2}_{xx} \varphi(\cdot)] +2q\lambda^{\top} \sigma^{\top} \partial^{2}_{xq}\varphi(\cdot) + |\lambda|^2q^2\partial^{2}_{qq}\varphi(\cdot)  \right)=0\,,\label{eq ViscSol}
       \end{gather}
  on the different time steps, and of the following boundary condition obtained by taking the Fenchel transform in  \eqref{eq: Hjb bord1}
  \begin{gather}
        v^\sharp(\tip-,\cdot)=(v\veee \ell)^{\sharp}(\tip,\cdot)\,. \label{eq ViscSol bord}
   \end{gather}
 By the Feynman-Kac representation this corresponds to the following representation
 \begin{eqnarray*}
 v^\sharp(t,\cdot)&\le&\E^{\Q_{t,x}}\left[(v\veee \ell)^{\sharp}(\tip,\cdot)\right]\,, \;\text{ for }\; t\in [t_{i},t_{i+1})\,,\;i<n\,.
 \end{eqnarray*}
 The aim of this paper is actually to prove by using probabilistic arguments only that on $\Rpp\times\R$
 \begin{eqnarray}\label{eq: def w}
 \left\{ \begin{array}{rcl}
  {w(T,x,q)}&:=& {q+\infty\1_{\{q<0\}}} \,,\\
 w(t,x,q)&:=&\E^{\Q_{t,x}}\left[(w^{\sharp}\veee \ell)^{\sharp}(\tip,X^{t,x}_{\tip},Q^{t,x,q}_{\tip})\right]\,, \;\text{ for }\; t\in [t_{i},t_{i+1})\,,\;i<n\,,
 \end{array}\right.
 \end{eqnarray}
with  $Q^{t,x,q}$ defined in \eqref{eq: def Q},
is the proper algorithm to compute the value function $v^\sharp$ and thus $v$.

 Indeed our main result is given by the following theorem. 
  \begin{theo}\label{th mainrsult} $v=w^{\sharp}$ on  {$[0,T]\x \Rpp\x [0,1]$}.
  \end{theo}

The proof of this result is the object of the subsequent sections.
Although it is in the spirit of \cite{bouchard2009stochastic}, our proof is
different and more involved. The main difficulty comes from the induction. At
each time step, we have to verify that $(w^{\sharp}\veee \ell)$ behaves in a
sufficiently nice way. In the one step case,  \cite{bouchard2009stochastic} had
only to consider the terminal payoff $\ell$. Moreover, we only use probabilistic
arguments as opposed to PDE arguments.

 \vspace{2mm}

{Clearly, the algorithm \eqref{eq: def w} provides a way to compute the value
function easily. One can for instance use the fact that $w=v^{\sharp}$ is the
unique viscosity solution of \eqref{eq ViscSol} with the boundary condition \eqref{eq ViscSol bord}. Let us make  this statement more precise.}
 \begin{dfn} We say that a lower-semicontinuous function $u$ is a viscosity super-solution of the system ($\cal S$) if, on each  $[\ti,\tip)\x\Rpp\x (0,\infty)$, $i<n$, it is a viscosity super-solution of
\eqref{eq ViscSol}
  with  the boundary conditions
 \begin{eqnarray*}
 \begin{array}{rcl}
      \liminf\limits_{t'\uparrow \ti, (x',q')\to (x,q)}  u(t',x',q')&\ge&
(u^{\sharp}\veee \ell)^{\sharp}(t_i,x,q)\;\text{ for }\; (x,q)\in \Rpp\x
(0,\infty)\,,\; i<n\,,\\
       \liminf\limits_{t'\uparrow T, (x',q')\to (x,q)}  u(t',x',q')&\ge&  \ell^{\sharp}(T,x,q)\;\text{ for }\; (x,q)\in \Rpp\x (0,\infty)\,.
       \end{array}
  \end{eqnarray*}
We define accordingly the notion of sub-solution for  upper-semicontinuous functions. A function is a viscosity solution if its lower- (resp. upper-) semicontinuous envelope is a viscosity super- (resp. sub-) solution.
 \end{dfn}
 Note that in the above definition we have to understand $u$ as being $+\infty$  on $[0,T]\x\Rpp\x(-\infty,0)$ to compute the Fenchel transforms involved in  the time boundary conditions.
\vspace{2mm}

{We now provide a version of the comparison principle for ($\cal S$) which
pertains for the usual extensions of the Black and Scholes model. The
assumptions used below are here to avoid the boundary of $\Rpp$ - when this is
not the case, one has to specify additional boundary conditions.}

 \begin{prop}\label{prop: solution visco et unicite} The function $w$ is continuous on {$[t_i,t_{i+1})\x \Rpp\x \R_{+},\,i<n$}, non-negative, has linear growth in its last variable and is a viscosity solution of ($\cal S$). Moreover, if    there exists two   functions $\bar \sigma$ and $\bar \mu$ such that $\sigma(\cdot,x)={\rm diag}[x]\bar \sigma(\cdot,x)$ and $\mu(\cdot,x)={\rm diag}[x]\bar\mu(\cdot,x)$, then  $u_{1}\ge u_{2}$ on  $[0,T)\x\Rpp\x {{(0,\infty)}}$  whenever  $u_{1}$ and $u_{2}$ are respectively a super- and a sub-solution of ($\cal S$), which are non-negative and have linear growth in their last variable on $[0,T)\x \Rpp\x \R_{+}$.  \end{prop}

%
%
%
%

The proof is postponed to the Appendix. Given the latter, it is not difficult to
follow the arguments of \cite{barles1991convergence} to construct a convergent
finite difference scheme for the resolution of $(\cal S)$. Alternatively, one could
also use quantization methods to tackle the approximation  of $w$, see
\cite{balpag03,balpag05}, or a regression based Monte-Carlo method, see the
survey paper \cite{bouchard2012monte} and the references therein.

\subsection{Examples of application}

In this section, we present two examples of application. The numerical results
are obtained using the following procedure which is based on the above
algorithm to compute $w=v^\sharp$: for $i \le n-1$,\\
1) compute the value of $(w^\sharp\vee \ell)^\sharp(t_{i+1},\cdot)$ by
approximating the Fenchel-Legendre transform numerically, \\
2) solve the PDE \eqref{eq ViscSol}-\eqref{eq ViscSol bord} for $w$,
using e.g. finite difference methods, on $[t_i,t_{i+1}]\times\Rpp\times\R_+$.\\
%

We now fix $T=1$ and $\bT_t\defeq \set{t_0=0,\,t_1=\frac{1}{3},\,t_2=\frac{2}{3},\,t_3=1}\cap(t,t_3],\,t\in[0,T]$.
We work in a Black-Scholes setting with market parameters: $d=1$,
$\sigma(t,x)=0.25x$, $\lambda(t,x)=0.2$.

\vspace{2mm}
For our first numerical application, we  consider a put option, i.e. $\ell(t,x)=[K-x]^{+}$, with strike $K=30$.

\begin{figure}[H]
     \centering
     \begin{tabular}{cccc}
      \includegraphics[width=.5\linewidth]{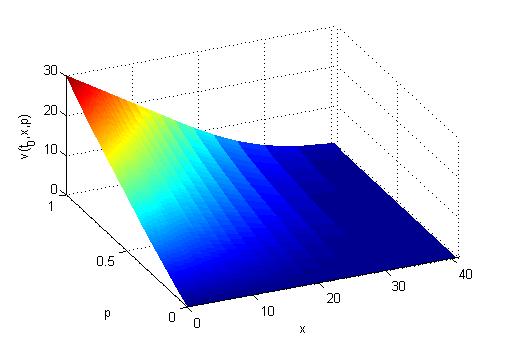} &
      \includegraphics[width=.5\linewidth]{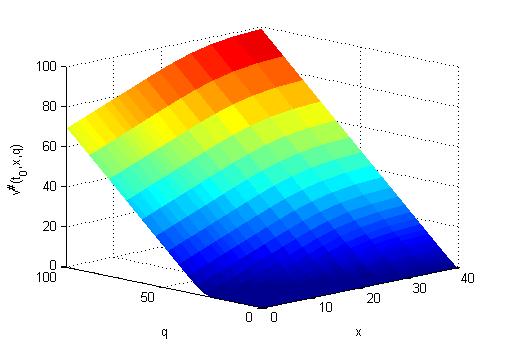}\\
      (a)&(b)
     \end{tabular}
     \caption{Surface of $v(t,x,p)$ and $v^\sharp(t,x,q)$ at $t=t_0$.\label{fig:put1}}
   \end{figure}

\begin{figure}[H]
      \centering
      \begin{tabular}{cccc}
       \includegraphics[width=.5\linewidth]{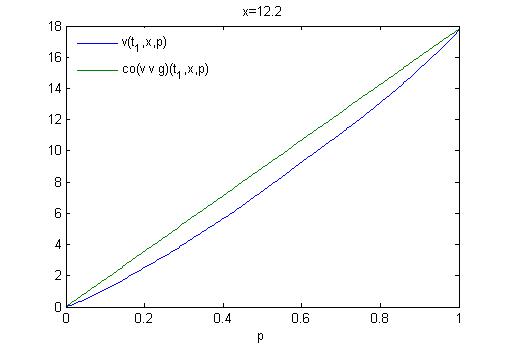} &
       \includegraphics[width=.5\linewidth]{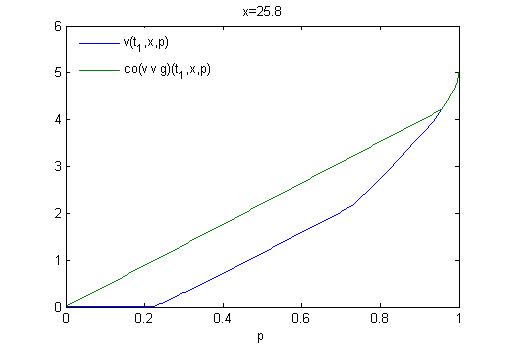} \\
       (a)&(b)\\
       \includegraphics[width=.5\linewidth]{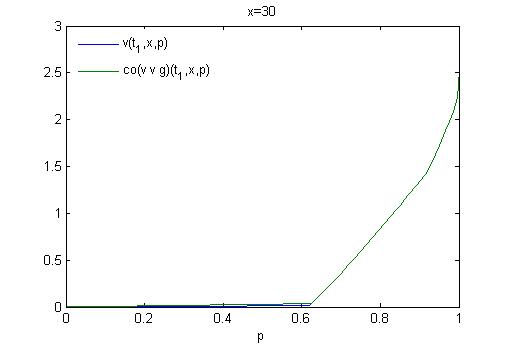} &
       \includegraphics[width=.5\linewidth]{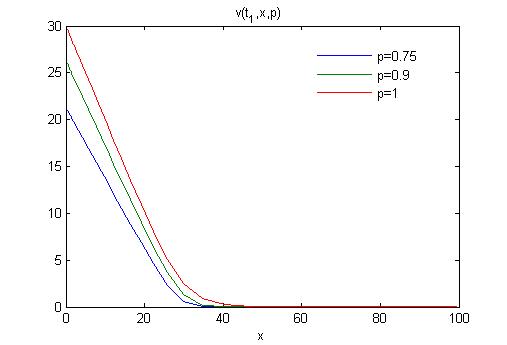} \\
       (c)&(d)
      \end{tabular}
      \caption{(a)-(c): plots of $v(t,x,\cdot)$ and $\co{v \vee
\ell}(t,x,\cdot)$ at $t=t_1$ and for different values of $x$. (d): plot
of $v(t,\cdot,p)$ at $t=t_1$ and for different values of $p$. \label{fig:put2}}
    \end{figure}

In figure \ref{fig:put1}, we plot  the functions $v$ and $v^\sharp$ at
$t=t_0$. In figure  \ref{fig:put2}(a-b-c), we plot for
different values of $x$ the function $v$ and $\co{v \vee \ell}$.
 This shows the rather complicated behavior  of the transformation $v \mapsto
\co{v\vee \ell}$, as predicted by Proposition \ref{thm
explicittildek}(b) {below}. With the notation of this proposition, figure
\ref{fig:put2}(a) corresponds to the case $A_1$, figure \ref{fig:put2}(b)
corresponds to the case $A_3$ and figure \ref{fig:put2}(c) corresponds to the
case $A_2$. Because of the interest rate being set to $0$ and the payoff being
convex, we always have $v(t,x,1) \ge \ell(t,x)$. Figure \ref{fig:put2}(d) shows
the decrease of value for $v$, when $p$ decreases.

\begin{figure}[H]
      \centering
      \begin{tabular}{cccc}
       \includegraphics[width=.5\linewidth]{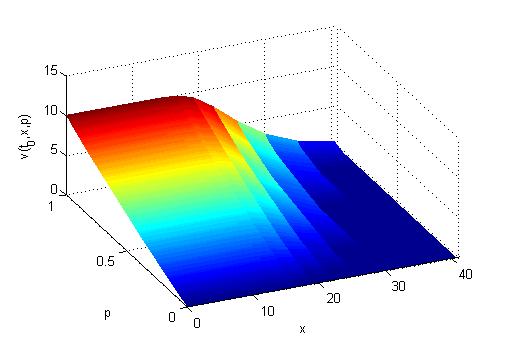} &
       \includegraphics[width=.5\linewidth]{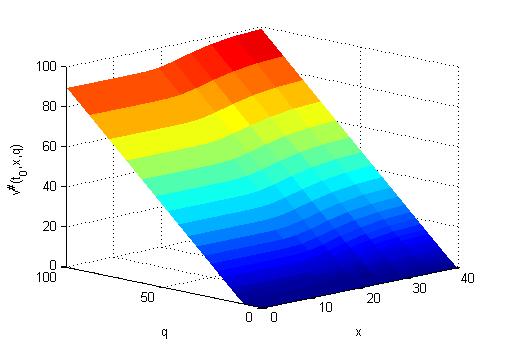}\\
      (a)&(b)
      \end{tabular}
      \caption{Surface of $v(t,x,p)$ and $v^\sharp(t,x,q)$ at $t=t_0$.\label{fig:spreadput1}}
   \end{figure}

\begin{figure}[H]
    \centering
    \begin{tabular}{cccc}
      \includegraphics[width=.5\linewidth]{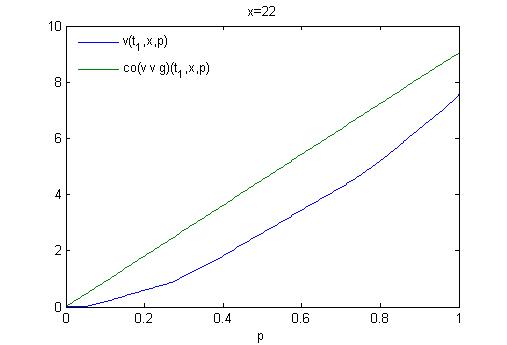} &
      \includegraphics[width=.5\linewidth]{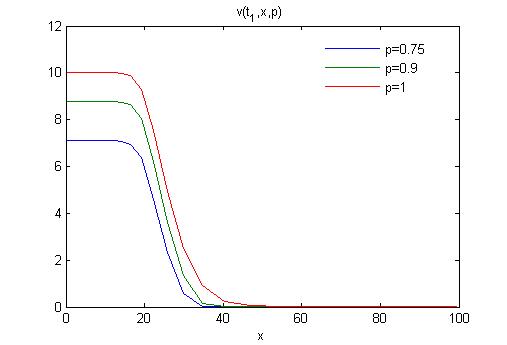}\\
      (a)&(b)
    \end{tabular}
    \caption{(a): plot of $v(t,x,\cdot)$ and $\co{v \veee g}(t,x,\cdot)$ at
$t=t_1$ and for $x=22$. (b): plot of $v(t,.,p)$ at $t=t_1$ and for different
values of $p$.\label{fig:spreadput2}}
  \end{figure}

 In our second example, we consider a put spread option with strikes $20$ and $30$, i.e. $\ell(t,x)=[30-x]^{+}-[20-x]^{+}$. The numerical results are displayed in Figure \ref{fig:spreadput1} and \ref{fig:spreadput2}. It may happen here that $v(t,x,1) < g(t,x)$, see figure \ref{fig:spreadput2}(a).

We conclude this section {with} the following remark on the behavior of $v$
near ${p=}1$.
\begin{rem}
(a) We know from the identification $v=w^\sharp$ and Proposition \ref{pr basic k and tilde k}(b) that $p\mapsto v(t,x,p)$ is convex and continuous on $[0,1]$.
\\
(b) Nothing prevents $\mathrm D^-_pv(\cdot,1)$ to be equal to $+\infty$. This can be checked by direct calculation in the European case and the Black-Scholes setting using the explicit formula  \cite[Equation (3.15)]{follmer1999quantile}.
\end{rem}

 \section{Proof of the backward  dual representation}\label{sec: proofs}

From now on, we extend $v$ to $[0,T]\x \Rpp\x\R$ by setting
\begin{equation}\label{eq: extension v a R}
v(\cdot,p)=0 \text{ if } p<0 \text{ and } v(\cdot,p)=+\infty \text{ if } p>1\,.
\end{equation}
Using the convention $\inf \emptyset =+\infty$, this extension is consistent with \eqref{eq: def v}.

 \subsection{The backward algorithm as a lower bound}

 We first show that the backward algorithm \eqref{eq: def w} actually provides a lower bound for the value function $v$.

       \begin{prop} \label{pr w ge ksharp} $v\ge w^{\sharp}$ on  {$[0,T]\x \Rpp\x [0,1]$}.
      \end{prop}

      \proof First note that  {$v(T,\cdot)=0=w^{\sharp}(T,\cdot)$}, by definition. Thus, $(v\vee \ell)(T,\cdot) =(w^{\sharp}\vee \ell)(T,\cdot)$. We now assume that  {$v\ge w^{\sharp}$} on $[\tip,T]\x \Rpp\x [0,1]$ for some $i\le n-1$.
      Then, $(v\veee \ell)^{\sharp}(\tip,\cdot)\le (w^{\sharp} \veee \ell)^{\sharp}(\tip,\cdot)$ and therefore
         \begin{align*}
         (v\veee\ell)\left(\tip,X^{t,x}_{\tip},P^{t,p,\alpha}_{\tip}\right) &
         \ge P^{t,p,\alpha}_{\tip} qQ^{t,x,{1}}_{\tip} -
 (v\veee\ell)^\sharp\left(\tip,X^{t,x}_{\tip},qQ^{t,x,{1}}_{\tip}\right)
         \\
         & \ge P^{t,p,\alpha}_{\tip} qQ^{t,x,{1}}_{\tip} -
 (w^{\sharp}\veee\ell)^{\sharp}\left(\tip,X^{t,x}_{\tip},qQ^{t,x,{1}}_{\tip}\right)\,.
        \end{align*}
    Fix $t\in [\ti,\tip)$.    Taking the expectation on both sides and recalling \eqref{eq: def w}, we obtain
        \begin{align*}
\E^{\Q_{{t,x}}}\left[(v\veee\ell)\left(\tip,X^{t,x}_{\tip},P^{t,p,\alpha}_{\tip}\right)\right] \ge pq-w(t,x,q)\,.
        \end{align*}  Taking first the supremum
over $q\in\R$ in the right-hand side  and then
        the infimum over $\alpha\in\mathcal{A}_{t,p}$ in the left-hand side, we
get from Proposition \ref{wProbRep0} that $v(t,x,p)\ge w^\sharp(t,x,p)$.      \eproof

 \subsection{Representation and differentiability of  the backward dual algorithm}\label{sec: travaille sur w}

 This section is devoted to the study of the function $(w^{\sharp}\veee
\ell)^{\sharp}$ which appears in the dual algorithm \eqref{eq: def w} and of its
Fenchel transform $(w^{\sharp}\veee \ell)^{\sharp\sharp}$.
We first provide a decomposition in simple terms in Proposition \ref{thm
explicittildek}. They only contain $w, \ell$ and auxiliary functions that are
easy to handle, see \eqref{eq ktildesharp0}-\eqref{eq ktilde0} below. In view of \eqref{eq: def w},
this will then allow us to study the subdifferential of $w(\ti,\cdot)$ in terms
of the subdifferential of $w(\tip,\cdot)$. This analysis is reported in Lemma
\ref{le prop Dk}. These results will be of important use in the final proof of
Theorem \ref{th mainrsult} as it will require to find a particular value $p$ in
the subdifferential of $w(\ti,\cdot)$ and then to apply a martingale
representation argument between elements of the subdifferential of
$(w^{\sharp}\veee \ell)^{\sharp}$ at $\tip$ and $p$ at $\ti$, see the proof of
Theorem \ref{kRepnew}.

\vspace{2mm}


 We start with properties that stem directly  from the definition of $w$  and
standard results in convex analysis. The proof is postponed to the Appendix.

    \begin{prop}\label{pr basic k and tilde k}
{The following holds for all $(t,x) \in
[0,T]\times\Rpp$.}

{\rm (a)} The functions $q\in \R  \mapsto w(t,x,q)$ is  a
proper convex non-decreasing and non-negative function. Moreover,
$w(\cdot,0)=0$ and $w(\cdot,q)=\infty$ for $q<0$.

{\rm (b)} The function $p\in \R  \mapsto w^\sharp(t,x,p)$ and   $q \in \R
\mapsto (w^{\sharp}\veee \ell)^{\sharp}(t,x,q)$  are   convex, non-negative,
non-decreasing and continuous on
their respective domains. Moreover, $w^\sharp(\cdot,0)=0=(w^{\sharp}\veee
\ell)^{\sharp}(\cdot,0)$ and $(w^{\sharp}\veee
\ell)^{\sharp}(\cdot,q)=+\infty$ for $q < 0$.
  \end{prop}
%



The next result is key to get the representation of $(w^\sharp\vee g)^\sharp$
and $(w^\sharp\vee g)^{\sharp \sharp}$.
 Recall that $\ell(t,x,p)=\ell(t,x)\1_{\{0<p\le 1\}}+\infty\1_{\{p> 1\}}$.
  \begin{lem} \label{rem: env conv de ell tronquee}
Let $p_1 \ge 0$ and $f$ be a non-decreasing convex function such that $f(0)=0$,
$f\ge \ell(t,x,\cdot)$ on $[p_{1},\infty)$, $f\le \ell(t,x,\cdot)$ on
$(-\infty,p_{1}]$.
\\
  {\rm (a)}  The convex envelope of $f\vee \ell$ is given by
  $$
  {(f\veee\ell)^{\sharp\sharp}(p)=}\co{f\veee\ell }(p)= pq_1\1_{\{0\le p< p_{1}\}}+ f(p)\1_{\{ p_{1} \le p \le
1\}} + \infty \1_{\{{p>1}\}}\,,$$
  with $q_1 = \ell(t,x)/p_1\1_{\set{p_1>0}}$.
\\
  {\rm (b)} Moreover, we have
  \begin{align*}
  (f\veee\ell )^\sharp(\cdot,q) = p_1[q-q_1]^{+}\1_{\set{q {\le} \mathrm D^+_pf(p_1)}} +
f^\sharp(q)\1_{\set{q{>} \mathrm D^+_pf(p_1)}}\,,\;q \ge 0\,,
  \end{align*}
  which is a closed proper convex function. In particular, it is continuous at
$\mathrm D^+_pf(p_1)$ when $0<\mathrm D^+_pf(p_1)<+\infty$.
  \end{lem}

  \proof \\
  \textbf{1.} {The left-hand side identity in (a) follows from \cite[Theorem
12.2]{rocka}.}
We {set} $\varphi:p\mapsto pq_1\1_{\{ p>0\}}\vee f(p)$, which is
convex.
 By assumption, we already know that $f(p) \le \ell(t,x,p)=0$ for $p \le 0$.
Since $f(0)=0$ and $f(p_1)=\ell(t,x)$, we have by convexity  that
$f(p) \le p q_1$, $p \in [0,p_1]$, which implies $\varphi(p)\1_{\set{p\le p_1}}
=pq_1\1_{\{0\le p\le p_{1}\}}$, for $p \le p_1$.
 Since $f(p)\le p q_1$ for $p \in [0,p_1]$ and $f(p_1)=p_1q_1$, we compute that
$\mathrm D^-_pf(p_1) \ge q_1$. By convexity, we also have
 $f(p) \ge f(p_1) +\mathrm D^-_pf(p_1)(p-p_1) \ge p q_1 $ for $p \ge p_1$ and then
$\varphi\1_{[p_1,\infty)} = f\1_{[p_1,\infty)}$. In particular, we observe that
$\varphi \le f\veee\ell $. It is straightforward to check that any candidate for
the convex envelope of $f\veee \ell$ is below $\varphi$. The above shows also that {$\mathrm D^+_pf(p_1)>
0$ whenever $q_{1}>0$.}
  \\
  \textbf{2.} Let us now observe that $f^\sharp(q)<\infty$, for $q\ge 0$, since $f(\cdot,p)= g(\cdot,p)=\infty$ for $p>1\vee p_1$. It follows that   the subdifferential of $f^\sharp$ at non-negative $q$ is non empty. The
proof of (b) follows from calculations based on the following results from
convex analysis, see e.g.~\cite[Chapter I Proposition 5.1]{eketem76}. Let
$\psi$ be a proper function on $\R$, then
  $p$ is in the subdifferential of $\psi$ at $q$ if and only if
\begin{align}\label{eq duality for dummies 1}
 \psi^\sharp(p) + \psi(q) = pq \,. 
  \end{align}
  (i) At $p=0$, the subdifferential of  $(f\veee\ell )^{\sharp \sharp} =
\co{f\veee\ell }$ is equal to $[0,q_1]$. This follows directly from the
characterization of the convex envelope of $f\veee\ell $ given in (a). Using
the above equality with $\psi = (f\veee\ell )^\sharp$, we then have for $q
\in  [0,q_1]$
  \begin{align*}
  (f\veee\ell )^{\sharp \sharp}(0) + (f\veee\ell )^\sharp(q) = 0\times q \implies
(f\veee\ell )^\sharp(q) = 0\,,
  \end{align*}
  {since $(f\veee\ell )^{\sharp \sharp}(0)=0$ by our assumption,
  namely $f(0)=0=g(\cdot,0)$ and $ g\ge 0$.}\\
  (ii)  {The subdifferential of $(f\veee\ell )^{\sharp \sharp}=\co{f\veee\ell }$ at
$p_1$  is equal to $\mathcal{D}:=[q_1,\mathrm D^+_pf(p_1)]$} if
  $\mathrm D^+_pf(p_1)<+\infty$ or $[q_1,+\infty)$ otherwise. This follows again
directly from (a). We recall from the step 1.~that $f(p_{1})=q_{1}p_{1}$. Then,
  using \eqref{eq duality for dummies 1} with $\psi = (f\veee\ell )^\sharp$  {and (a)}, we
 have for   $q \in \mathcal{D}$
    \begin{align*}
  (f\veee\ell )^{\sharp \sharp}(p_1) + (f\veee\ell )^\sharp(q) = p_1  q \implies
(f\veee\ell )^\sharp(q){=p_{1}q-f(p_{1})} = p_1  (q -  q_1){=p_1  [q -
q_1]^{+}}\,.
  \end{align*}
  (iii) If {$ q > \mathrm D^+_pf(p_1)$}, {an element $p$ of the  subdifferential of $f^\sharp$ at $q$ satisfies}
    \begin{align*}
 f(p) + f^\sharp(q) = p q \,.
  \end{align*}
 We   first note that  $p \ge p_1$ necessarily. Indeed, by  \cite[Chapter I
Corollary 5.2]{eketem76}, $q\in [\mathrm D^{-}_{p}f(p),\mathrm D^{+}_{p}f(p)]$ while $q > \mathrm D^+_pf(p_1)$. Recall that
$f=(f\veee\ell )^{\sharp \sharp}$ on $[p_1,\infty)$. We then deduce from the
previous
equality  that
  \begin{align*}
  (f\veee\ell )^{\sharp \sharp}(p) +f^\sharp(q) = p q  \implies  f^\sharp(q){=pq-(f\veee\ell )^{\sharp \sharp}(p)} \le
(f\veee\ell )^\sharp(q)  \,.
  \end{align*}
  Observing that the reverse inequality follows from $f \le f\veee\ell $, we
get $ f^\sharp(q) = (f\veee\ell )^\sharp(q)$ for $q \in  {(\mathrm D^+_pf(p_1)},+\infty)$.
  \eproof

  \vspace{2mm}

We are now in position to provide the decomposition of
$(w^\sharp\vee g)^\sharp$ and
$(w^\sharp\vee g)^{\sharp\sharp}$. It basically follows from the application of
the previous Lemma to  $f=w^\sharp$.
  \begin{prop}\label{thm explicittildek} For $(t,x,p) \in
[0,T]\times\Rpp\times\R$, we define the following `facelift' of $\ell$
 \begin{align*}
     \zeta(t,x,p) = q_\ell(t,x) p\1_{\{0 \le p \le 1 \}}
+\infty\1_{\{p > 1\}}  \,.
    \end{align*}
    with
    \begin{align*}
    q_\ell(t,x) := \frac{\ell(t,x)}{p_\ell(t,x)} \1_{\set{p_\ell(t,x)>0}} \;
\mbox{ and }\;
      p_\ell(t,x) := \sup \set{p \in \R \,|\, w^\sharp(t,x,p)=\ell(t,x) }\;
\wedge 1\,.
    \end{align*}
    Then,
      \begin{enumerate}[{\rm (a)},noitemsep,nolistsep]
        \item  The function $q \mapsto (w^{\sharp}\veee
\ell)^{\sharp\sharp}(\cdot,q)$ is continuous on its domain and
         \begin{align}\label{eq ktildesharp0}
          (w^{\sharp}\veee \ell)^{\sharp\sharp} & =
\co{w^{\sharp}\veee \ell}
           = w^{\sharp}\veee \zeta\,.
         \end{align}
         \item For all $q\in \R_{+}$:
         \begin{align}
         (w^{\sharp}\veee \ell)^{\sharp}(\cdot,q)& =
\left[q-\ell(\cdot)\right]^+ \1_{A_1}(\cdot)
         + w(\cdot,q) \1_{A_2}(\cdot)+
          \kappa(\cdot,q)\1_{A_3}(\cdot)\,,\label{eq ktilde0}
         \end{align}
         where {
         \begin{align*}
         \kappa(\cdot,q) &:=
p_\ell(\cdot)\left[q-q_\ell(\cdot)\right]^+\1_{\set{q\le\bar{q}(\cdot)}}
         + w(\cdot,q) \1_{\set{q> \bar{q}(\cdot)}}\,,
        \end{align*}
with
$\;\bar{q}(\cdot):=\mathrm D^+_pw^\sharp\left(\cdot,p_\ell(\cdot)\right)$ and the subsets of $\,[0,T]\times \Rpp$:
$A_1 =\set{\ell > 0,\; {w^\sharp(\cdot,1)\le
\ell} }$, $A_2=\set{ \ell = 0}$, $A_3=\set{\ell> 0, \; {w^\sharp(\cdot,1) >
\ell}}$.
}

        \end{enumerate}
\end{prop}

    \begin{rem}\label{re p}
{\rm (a)}   It follows from Proposition \ref{pr basic k and tilde k}   that
$w^{\sharp}(\cdot,0)=0$. Hence,  $\ell(t,x)>0$ implies
$p_\ell(t,x)>0$ and
        \begin{align*}
         q_\ell(t,x)& =
\frac{\ell(t,x)}{p_\ell(t,x)}\1_{\set{\ell(t,x)>0}} \text{ so that }
 q_\ell(t,x)=0 \mbox{ if and only if } \ell(t,x)=0 \,.
        \end{align*}
        {\rm (b)} The decomposition on $A_1$, $A_2$ and $A_3$ will be
useful in the sequel, see e.g. proof of Lemma \ref{le prop Dk}(c) below.
\\
{\rm (c)} On $A_3$, we have $\bar{q} > 0$  {since $w^{\sharp}(\cdot,
p_{g}(\cdot))\ge g> 0$ and $w^{\sharp}(\cdot,0)=0$,
see Proposition \ref{pr basic k and tilde k}.}

    \end{rem}

{\bf Proof of Proposition \ref{thm explicittildek}.}
  {The identities in  \eqref{eq ktildesharp0} are   immediate consequences} of Lemma \ref{rem: env conv de ell tronquee}(a), Proposition \ref{pr basic k and tilde k}(b)  and
of the definition of $p_g$.
We now prove \eqref{eq ktilde0}. For $(t,x)\in A_1$, we {have $w^{\sharp}(t,x,\cdot)\le g$ and therefore $(w^\sharp\veee \ell)^{\sharp}(t,x,\cdot) = \ell^{\sharp}(t,x,\cdot)$ $= $ $[\cdot-\ell(t,x)]^{+}$ on $\R_{+}$.}
%
For $(t,x) \in A_2$, we have that $w^\sharp \ge g$ by Proposition
\ref{pr basic k and tilde k}(b) and the result follows directly. On
$A_3$, the expression is exactly the one given by Lemma  \ref{rem:
env conv de ell tronquee}(b).
\eproof

\vspace{4mm}

We can now turn to the study of the subdifferential of $w$. Recall the
definition of $p_{\rm min}$  in \eqref{eq: def pmin}.
   \begin{lem}\label{le prop Dk}
      Fix $0\le i\le n-1$ and $(t,x)\in[t_i,t_{i+1})\times\Rpp$. Then:\\
         {\rm (a)} $\mathrm D^+_q w(t,x,\cdot)\ge 0$ if $q\ge 0$ and
 $\mathrm D^-_q w(t,x,\cdot)\ge 0$ if $q>0$,\\
        {\rm (b)} {$\lim_{q\uparrow\infty}\mathrm D^{+}_q w(t,x,q)= 1$},\\
        {\rm (c)} $\mathrm D^+_q w(t,x,0)=p_{\rm min}(t,x)$.\\
      Moreover,
      \begin{align}
 \mathrm D^-_qw(t,x,q)&=
\E\left[\mathrm D^-_{q}(w^\sharp\veee
\ell)^\sharp(t_{i+1},X^{t,x}_{\tip},qQ^{t,x,1}_{\tip}
))\right]\,\; \text{ for } q > 0\,, \; \text{ and } \label{eq express D-k}
\\
          \mathrm D^+_qw(t,x,q)&=
\E\left[\mathrm D^+_{q}(w^\sharp\veee
\ell)^\sharp(t_{i+1},X^{t,x}_{\tip},qQ^{t,x,1}_{\tip}
))\right]\,\; \text{ for } q \ge 0\,.\label{eq express D+k}
\end{align}
    \end{lem}

    \proof The proof is based on an induction argument. Our assumptions
guarantee that (a)-(b)-(c) are valid at $T$. Let us assume that it holds true on
$[\tip,T]$ for some $i\le n-1$.

In view of  Proposition \ref{thm explicittildek}, we obtain for $q\ge 0$ and
$j\le n$ that
        \begin{align*}
         \mathrm D^+_q(w^\sharp\veee \ell)^\sharp(t_{j},x,q)& =
\1_{\set{q \ge \ell(t_{j},x)}} \1_{A_1}(t_{j},x)
 +
\mathrm D^+_q w(t_{j},x,q) \1_{A_2}(t_{j},x)
\\&\quad+
\mathrm D^+_q\kappa(t_{j},x,q) \1_{A_3}(t_{j},x)\,,
        \end{align*}
with
\begin{align*}
 \mathrm D^+_{q}\kappa(t_{j},x,q) = p_\ell(t_{j},x)\1_{\set{q_\ell(t_{j},x) \le q <
\bar{q}(t_{j},x)}} + \mathrm D^+_{q}w(t_{j},x,q)\1_{\set{q >
\bar{q}(t_{j},x)}}\,.
\end{align*}
For $q >0$, we have
        \begin{align*}
         \mathrm D^-_q(w^\sharp\veee \ell)^\sharp(t_{j},x,q)& =
	\1_{\set{q > \ell(t_{j},x)}} \1_{A_1}(t_{j},x)
         +
         \mathrm D^-_q w(t_{j},x,q) \1_{A_2}(t_{j},x)
\\&\quad+
\mathrm D^-_q\kappa(t_{j},x,q) \1_{A_3}(t_{j},x)\,,
        \end{align*}
with
\begin{align*}
 \mathrm D^-_{q}\kappa(t_{j},x,q) = p_\ell(t_{j},x)\1_{\set{q_\ell(t_{j},x) < q \le
\bar{q}(t_{j},x)}} + \mathrm D^-_{q}w(t_{j},x,q)\1_{\set{q >\bar{q}(t_{j},x)}}\,.
\end{align*}

We have by induction $\lim_{q \uparrow \infty} \mathrm D^+_{q}
\kappa(\tip,x,q) = 1$, which ensures that $\lim_{q \uparrow \infty}
\mathrm D^+_q(w^\sharp\veee \ell)^\sharp(\tip,x,q) = 1$. By the convexity of
$(w^\sharp\veee \ell)^\sharp$, this implies that $\mathrm D^+_q(w^\sharp\veee
\ell)^\sharp(\tip,x,q) \le 1$.
In view of \eqref{eq: def w}, a dominated convergence argument then leads to
\eqref{eq express D-k}-\eqref{eq express D+k} and   $\lim_{q \uparrow +\infty}
\mathrm D^+_{q}w(t,x,q) = 1$.

We now use our
induction hypothesis again to observe from the decomposition above that
\begin{align*}
\mathrm D^-_q(w^\sharp\veee \ell)^\sharp(\tip,x,q) \ge 0 \, ,\; q>0\,, \; \text{ and }\;
\mathrm D^+_q(w^\sharp\veee \ell)^\sharp(\tip,x,q) \ge 0 \,,\; q \ge 0\,.
\end{align*}
Recalling \eqref{eq express D-k}-\eqref{eq express D+k}, this shows  that
$\mathrm D^-_qw(t,x,q) \ge 0$ for $q>0$ and  $ \mathrm D^+_qw(t,x,q) \ge 0$ for $q \ge 0$.

It remains to prove (c).
From Remark \ref{re p}(a) and (c), the above decomposition implies  that
        $ \mathrm D^+_q(w^\sharp\veee \ell)^\sharp(\tip,x,0)$ $=$ $\mathrm D^+_q
w(\tip,x,0)\1_{\{\ell(\tip,x)=0\}}$. By our induction hypothesis, the last term
is  $ \mathrm D^+_q(w^\sharp\veee \ell)^\sharp(\tip,x,0)$ $=$ $p_{\rm
min}(\tip,x)\1_{\{\ell(\tip,x)=0\}}$. This identity combined with
\eqref{eq express D+k} provides
        \begin{align*}
 \mathrm D^+_qw(t,x,0) &=\esp{p_{\rm
min}(\tip,X^{t,x}_{\tip})\1_{\set{\ell(\tip,X^{t,x}_{\tip})=0}}}=p_{\rm
min}(t,x)\,,
\end{align*}
in which the last identity is an obvious consequence of the definition of
$p_{\rm min}$ in \eqref{eq: def pmin}.     \eproof

\begin{rem}\label{rem : sous gradient w en 0} Note that the subdifferential of
$w(t,x,\cdot)$ at $0$ is $(-\infty,p_{{\rm min}}(t,x)]$, since
$w(t,x,q)=\infty$ for $q<0$ and $\mathrm D^{+}_{q} w(t,x,0)=p_{{\rm min}}(t,x)$. See (a)
of Proposition \ref{pr basic k and tilde k}
and (c) of Lemma \ref{le prop Dk}.
\end{rem}

 \subsection{The backward algorithm as an upper-bound}

 Our final proof will  proceed by backward induction on the time steps. Fix $0\le i\le n-1$.
 Part of the induction hypothesis is:
 \begin{hypot}[\Hip] The following holds
  \begin{enumerate}[(i)]
\item\label{item : Hip co v sup ell continu} The functions $v(\tip,\cdot)$ and
$\co{v\veee \ell}(\tip,\cdot)$ are continuous on $\Rpp\x[0,1]$.
  \item\label{item : Hip co v sup ell bord 0 et 1}  $\co{v\veee
\ell}(\tip,\cdot,0)=0$ and $\co{v\veee \ell}(\tip,\cdot,1)=(v\veee
\ell)(\tip,\cdot,1)$.
 \item\label{item : Hip q-w} For all $x\in \Rpp$, the map $q\in \R_{+}\mapsto
q-(w^{\sharp}\veee \ell)^{\sharp}(\tip,x,q)$ is
{non-decreasing}, continuous and converges to $(v\veee
\ell)(\tip,x,1)$ as $q\to \infty$.
\end{enumerate}
 \end{hypot}

Before turning to the final argument, we provide three additional results that
hold at any time $ t \in [\ti{},t_{i+1})$ whenever \Hip~is in force.

 \subsubsection{Bounds and limits for $w^{\sharp}$}

 Our first additional result concerns the behavior of $w^{\sharp}$. It shows that $w^\sharp(\ti,x,1) = v(\ti,x,1)$. The last assertion will be used in
the proof of Lemma \ref{lem: induction Hip q-w} below to show that (iii) of
\Hi~holds if (iii) of \Hip~does.

 \begin{lem}\label{le more on ksharp}   Let {{\rm\ref{item : Hip q-w}}} of
\Hip~hold.
 Fix $(t,x) \in [\ti{},\tip)\times \Rpp$. Then,
 $w^{\sharp}(t,x,\cdot)$ is non-negative, continuous on its domain
$(-\infty,1]$ and
\begin{align*}
 0 \le w^\sharp(t,x,\cdot) \le w^\sharp(t,x,1) = v(t,x,1) \;\mbox{ on }
(-\infty,1]\,.
\end{align*}
Moreover,
the map $q\in \R\mapsto q-w(t,x,q)$ is non-decreasing, continuous on
$\R_{+}$ and converges to $v(t,x,1)$ as $q\to \infty$.
\end{lem}

\proof
The continuity and non-negativity of  $w^{\sharp}(t,x,\cdot)$ are stated in (b)
of Proposition \ref{pr basic k and tilde k}.
We now observe that \eqref{eq: def w} implies that
\begin{align*}
\delta(q):=q-w(t,x,q) = \E^{\Q_{t,x}}\left[ qQ^{t,x,1}_{\tip}- (w^{\sharp}\veee
\ell)^{\sharp} (\tip,X^{t,x}_{\tip},qQ^{t,x,1}_{\tip})\right]\,,
\end{align*}
which shows that $q\mapsto \delta(q)$ is non-decreasing since (iii) of \Hip~holds.
Applying the monotone convergence Theorem,
\ref{item : Hip q-w} of \Hip~and \eqref{eq: dpp retro pour v(t,x,1)}, we obtain
that $q\in \R_{+}\mapsto q-w(t,x,q)$ is   continuous and that
\begin{align*}
 \lim_{q \rightarrow \infty}\delta(q) = \E^{\Q_{t,x}}\left[
(v\veee \ell)(\tip,X^{t,x}_{\tip},1)\right] = v(t,x,1)\,.
\end{align*}
This implies that $w^\sharp(\ti,x,1) = \sup_{q\ge 0} \delta(q)$ $\ge$ $\lim_{q\to
\infty} \delta(q)$ $ =$ $ v(t,x,1)$, while  $w^\sharp(t,x,p)\ge  \lim_{q\to
\infty} (q(p-1)+\delta(q))$ $=$ $\infty$ for $p>1$. The fact that $w^\sharp(\ti,x,1) \le v(t,x,1)$ has been proved in Proposition \ref{pr w ge ksharp}.
\eproof

 \subsubsection{Convexification in the dynamic programming algorithms}

 As already mentioned in Remark \ref{rem convexifcation primal}, one can expect that $v\vee \ell$ can be replaced by its convex envelope, with respect to $p$, in \eqref{eq probw0}. The Hypotheses \ref{item : Hip co v sup ell continu}-\ref{item : Hip co v sup ell bord 0 et 1} of  \Hip~ensure this, see Proposition \ref{wProbRep1} below. We shall prove a similar result for $w^{\sharp}$ later on in Theorem \ref{kRepnew}. Note that the two identities \eqref{eq: v = inf E co p v sup ell} and \eqref{ProbRepksharp} below already suggest that the equality $v=w^{\sharp}$ at $\tip$ should iterate at $\ti$, since we already know from Proposition \ref{pr w ge ksharp} that $v\ge w^{\sharp}$.

    \begin{prop}\label{wProbRep1}   Let {\rm (i)-(ii)} of \Hip~hold. Then, for
all $t\in[t_i,t_{i+1})$ and  $(x,p)\in\Rpp\times[0,1]$, we have
      \begin{gather}
v(t,x,p)=\inf_{\alpha\in\mathcal{A}_{t,p}}\E^{\Q_{t,x}}\left[
\co{v\veee \ell}\left(t_{i+1},X^{t,x}_{t_{i+1}},P^{t,p,\alpha}_{t_{i+1}}\right)
        \right]\,.\label{eq: v = inf E co p v sup ell}
      \end{gather}
 Moreover, \bru{\rm\ref{item : Hip co v sup ell bord 0 et 1}} of \Hi~holds.
      \end{prop}

    \proof We fix $(t,x)\in[t_i,t_{i+1})\x\Rpp$.  Assuming that \eqref{eq: v = inf E co p v sup ell} is true, we   deduce that \ref{item : Hip co v sup ell bord 0 et 1} of \Hi~holds, since  $\Ac_{t,p}=\{0\}$  for $p\in \{0,1\}$ and therefore $P^{t,p,\alpha}_{\tip}=p$ for $\alpha\in \Ac_{t,p}$.  By \ref{item : Hip co v sup ell bord 0 et 1} of \Hip, the same argument combined with Proposition \ref{wProbRep0}  implies that \eqref{eq: v = inf E co p v sup ell} is valid for  $p\in \{0,1\}$.

    It remains to prove \eqref{eq: v = inf E co p v sup ell} for $0<p<1$. In view of  Proposition \ref{wProbRep0}, this reduces to showing that
\begin{eqnarray*}
  \inf_{\alpha\in\mathcal{A}_{t,p}}\E^{\Q_{t,x}}\left[
\co{v\veee \ell}\left(t_{i+1},X^{t,x}_{\tip},P^{t,p,\alpha}_{\tip}\right)
        \right]\ge \inf_{\alpha\in\mathcal{A}_{t,p}}\E^{\Q_{t,x}}\left[
 ({v\veee \ell})\left(t_{i+1},X^{t,x}_{\tip},P^{t,p,\alpha}_{\tip}\right)
        \right]\,,
\end{eqnarray*}
 the reverse inequality being trivial. We argue as in \cite[Proof of Proposition 3.3]{reveillac2012bsdes}.
It follows from the Caratheodory theorem  that we can find    two   maps $(\lambda_{j},\pi_{j}): (x,p)\in \Rpp\x [0,1] \mapsto  (\lambda_{j},\pi_{j})(x,p)\in  \Rpp\x [0,1]$, $j\le 2$, such that
\begin{eqnarray}\label{eq: point env convexe}
\begin{array}{c}
 \sum_{j=1}^{2} \pi_{j}(x,p)=1\mbox{ , }  p=\sum_{j=1}^{2}   \pi_{j}(x,p)\lambda_{j}(x,p) \\
  \mbox{ and }\; \co{v\veee \ell}(t_{i+1},x,p)=\sum_{j=1}^{2}   \pi_{j}(x,p)(v\veee \ell)(t_{i+1},x,\lambda_{j}(x,p))\,.
\end{array}
\end{eqnarray}
We claim that they can be chosen in a measurable way.
More precisely, \ref{item : Hip co v sup ell continu} of \Hip~and
\cite[Proposition 7.49]{bertsekas1978stochastic} imply that they can be chosen
to be analytically measurable. We can then appeal to \cite[Lemma
7.27]{bertsekas1978stochastic} to obtain a Borel-measurable version which
coincides a.e.~for the pull-back measure of
$(X^{t,x}_{\tip-\eps},P^{t,p,\alpha}_{\tip-\eps})$, for $\alpha\in \Ac_{t,p}$
and  $0<\eps<\tip-t$ fixed. This is this version that we use in the following.

 We now let $\xi$ be a $\Fc_{\tip}$-measurable random variable such that
$$
\P[\xi=\lambda_{j}(X^{t,x}_{\tip-\eps},P^{t,p,\alpha}_{\tip-\eps})|\Fc_{\tip-\eps}]=\pi_{j}(X^{t,x}_{\tip-\eps},P^{t,p,\alpha}_{\tip-\eps})\,.
$$
Then, $\E[\xi|\Fc_{\tip-\eps}]=P^{t,p,\alpha}_{\tip-\eps}$ by the above construction, and we can then find $\alpha_{\eps}\in \Ac_{t,p}$ such that
$P^{t,p,\alpha_{\eps}}_{\tip-\eps}=P^{t,p,\alpha}_{\tip-\eps}$ and $P^{t,p,\alpha_{\eps}}_{\tip}=\xi$. Recalling  \eqref{eq: point env convexe}, we obtain
\begin{align*}
  &\E\left[\left(Q^{t,x,1}_{t_{i+1}-\varepsilon}\right)^{-1}
\co{v\veee \ell}\left(t_{i+1},X^{t,x}_{\tip-\eps},P^{t,p,\alpha}_{\tip-\eps}\right)
        \right] \\&= \E^{\Q_{t,x}}\left[
 ({v\veee \ell})\left(t_{i+1},X^{t,x}_{\tip-\eps},P^{t,p,\alpha_{\eps}}_{\tip}\right)
        \right]\\&\quad- \E\left[\left(\left(Q^{t,x,1}_{t_{i+1}}\right)^{-1}-\left(Q^{t,x,1}_{t_{i+1}-\varepsilon}\right)^{-1}\right)({v\veee \ell})\left(t_{i+1},X^{t,x}_{\tip-\eps},P^{t,p,\alpha_{\eps}}_{\tip}\right)
        \right]\\
        &\ge
              \inf_{\alpha'\in \Ac_{t,p}}\E^{\Q_{t,x}}\left[
  ({v\veee \ell})\left(t_{i+1},X^{t,x}_{\tip},P^{t,p,\alpha'}_{\tip}\right)
         \right]+\Delta(\epsilon)\,,
\end{align*}
 with
 \begin{align*}
 \Delta(\epsilon)&=-C\E^{\Q_{t,x}}\left[
\left (1+|X^{t,x}_{\tip-\eps}|+|X^{t,x}_{\tip}|\right)\left|X^{t,x}_{\tip-\eps}-X^{t,x}_{
 \tip}\right|
   \right]\\&\quad-C\E\left[\left(\left(Q^{t,x,1}_{t_{i+1}}\right)^{-1}-\left(Q^{t,x,1}_{t_{i+1}-\varepsilon}\right)^{-1}\right)
  \left (1+|X^{t,x}_{\tip-\eps}|\right)\right]\,,
   \end{align*}
    recall \eqref{eq: ass ell lipschitz}, \eqref{eq: borne v} and \eqref{eq: v loc lip}.
Moreover, since $0\le \co{v \veee \ell}(t_{i+1},x,\cdot) \le {v \veee
\ell}(t_{i+1},x,\cdot)\le C(1+|x|)$, using \ref{item : Hip co v sup ell
continu} of  \Hip, we can pass to the limit  to obtain
\begin{eqnarray*}
  \E^{\Q_{t,x}}\left[
\co{v\veee \ell}\left(t_{i+1},X^{t,x}_{\tip},P^{t,p,\alpha}_{\tip}\right)
        \right]
        &\ge & \inf_{\alpha'\in \Ac_{t,p}}\E^{\Q_{t,x}}\left[
 ({v\veee \ell})\left(t_{i+1},X^{t,x}_{\tip},P^{t,p,\alpha'}_{\tip}\right)
        \right]\,.
\end{eqnarray*}
   \eproof

  Since our final result is $v=w^{\sharp}$, the same convexification should
appear in the dual algorithm.  As already mentioned, it will actually allow us
to show that $v=w^{\sharp}$ at $\ti$ if this true at $\tip$.
   \begin{theo}\label{kRepnew} Let {{\rm\ref{item : Hip q-w}}} of
\Hip~hold.
     Fix $(t,x,p)\in[\ti{},\tip)\times\Rpp\x[0,1]$. Then,
there exists  $\bar{\alpha}\in\mathcal{A}_{t,p}$ such that
      \begin{align}\label{ProbRepksharp}
       w^\sharp(t,x,p) &=
\E^{\Q_{t,x}}\left[\co{w^{\sharp}\veee \ell}\left(t_{i+1},X^{t,x}_{\tip},
P^{t,p,\bar{\alpha}}_{\tip} \right)\right]\,.
      \end{align}
          \end{theo}

    \proof Recall the definition of $p_{\rm min}$ in \eqref{eq: def pmin}.

    \textbf{1.} We first assume that  $p \in
(p_{\rm min}(t,x),1)$.
     We know from Lemma \ref{le prop Dk}(b)-(c) that there exists a
$\tilde{q}\in(0,\infty)$ such that $p$ lies in the subdifferential of
$w(t,x,\cdot)$ at $\tilde{q}$. Then, we can find
$\lambda\in[0,1]$ such that $p=\lambda \mathrm D^+_q
w(t,x,\tilde{q})+(1-\lambda)\mathrm D^-_q w(t,x,\tilde{q})$. In view of
\eqref{eq express D-k}-\eqref{eq express D+k}, this implies that
 \begin{align}
 p &= \E\left[(\lambda
\mathrm D^+_q(w^\sharp\veee \ell)^\sharp+(1-\lambda)\mathrm D^-_q(w^\sharp\veee
\ell)^\sharp)\left(t_{i+1},X^{t,x}_{\tip},\tilde{q}
Q^{t,x,{1}}_{\tip}\right)\right] \label{eq ineq1}\,.
     \end{align}
    It follows from Lemma \ref{le prop Dk} and its proof that the random
variable in the expectation is valued in $[0,1]$.
  By the martingale representation
theorem, we can find    $\bar{\alpha}\in \Ac_{t,p}$ such that
    \begin{gather*}
   (\lambda \mathrm D^+_q(w^\sharp\veee \ell)^\sharp+(1-\lambda)\mathrm D^-_q(w^\sharp\veee
\ell)^\sharp)\left(t_{i+1},X^{t,x}_{\tip},\tilde{q}
Q^{t,x,{1}}_{\tip}\right)
= p + \int_t^{t_{i+1}} \bar{\alpha}_{s}^{\top} \mathrm{d}
W_s =: P^{t,p,\bar{\alpha}}_{\tip}\,.
     \end{gather*}
 For later use, note that the above implies
     \begin{gather}\label{Probp}
P^{t,p,\bar{\alpha}}_{\tip} \tilde{q}Q^{t,x,{1}}_{\tip}
- (w^{\sharp}\veee \ell)^{\sharp}\left(t_{i+1},X^{t,x}_{\tip},\tilde{q}
Q^{t,x,{1}}_{\tip}\right)=(w^{\sharp}\veee
\ell)^{\sharp\sharp}\left(t_{i+1},X^{t,x}_{\tip},P^{t,p,\bar{\alpha}}_{\tip}
\right)\,,
 \end{gather}
where we used \eqref{eq duality for dummies 1} with $\psi = (w^\sharp\veee
g)^\sharp$.
On the other hand,   we also have, again by \eqref{eq duality for dummies 1}
with $\psi=w$,
     \begin{align}
      w(t,x,\tilde{q}) + w^\sharp(t,x,p)=\tilde{q}p\,,  \label{eq tfl prop}
     \end{align}
     and, by \eqref{eq: def w},
     \begin{gather}
      w(t,x,\tilde q) =
\E^{\Q_{t,x}}\left[(w^{\sharp}\veee
\ell)^{\sharp}\left(t_{i+1},X^{t,x}_{\tip},\tilde{q}
Q^{t,x,{1}}_{\tip}\right)\right]\,.\label{eq ineq0}
     \end{gather}

     Thus, inserting \eqref{eq ineq1}  and \eqref{eq ineq0}
into \eqref{eq tfl prop}, and using \eqref{Probp}, leads to
     \begin{align*}
   w^\sharp(t,x,p) &
   =
\E^{\Q_{t,x}}\left[
P^{t,p,\bar{\alpha}}_{\tip} \tilde{q}Q^{t,x,{1}}_{\tip}
- (w^{\sharp}\veee \ell)^{\sharp}\left(t_{i+1},X^{t,x}_{\tip},\tilde{q}
Q^{t,x,{1}}_{\tip}\right)\right]
   \\
&=\E^{\Q_{t,x}}\left[(w^{\sharp}\veee
\ell)^{\sharp\sharp}\left(t_{i+1},X^{t,x}_{\tip},
P^{t,p,\bar{\alpha } }_{\tip}\right)\right]\,.
     \end{align*}
 We conclude by appealing to \eqref{eq ktildesharp0}.

\textbf{  2.} We now assume that $p\in[0,p_{{\rm min}}(t,x)]$.
Since $[0,p_{{\rm min}}(t,x)] $ belongs to the subdifferential  of
$w(t,x,\cdot)$ at $0$, recall Remark
\ref{rem : sous gradient w en 0}, and    $p_{{\rm min}}(t,x)=\mathrm D^+_q
w(t,x,0)$, recall Lemma \ref{le prop Dk}, we can  find $\lambda \in [0,1]$ such
that
$ p = \lambda \mathrm D^+_q w(t,x,0)$. We then proceed as above up to obvious
modifications.

 \textbf{ 3.} We finally assume that $p=1$.
     We know from Lemma \ref{le more on ksharp} that $w^\sharp(t,x,1)=v(t,x,1)$.
Hence, \eqref{eq: dpp retro pour v(t,x,1)} implies
     \begin{align*}
       w^\sharp(t,x,1)=v(t,x,1) &= \E^{\Q_{t,x}}\left[(v\veee
\ell)\left(t_{i+1},X^{t,x}_{\tip},1\right) \right]\,.
          \end{align*}
As in the proof of Lemma \ref{le more on ksharp}, we deduce from \ref{item : Hip q-w} of \Hip that
  $\co{w^\sharp\veee \ell}(t_{i+1},\cdot,1) =(w^\sharp\veee \ell)^{\sharp\sharp}(t_{i+1},\cdot,1) \ge (v\veee
\ell)(t_{i+1},\cdot,1)$. In view of Proposition \ref{pr w ge ksharp}, this leads to   $ (v\veee \ell)(\tip,x,1)=\co{w^\sharp\veee \ell}(\tip,x,1)$.
\eproof

 \subsection{Conclusion of the proof}
    \label{subse end of proof}
 To conclude the proof of Theorem \ref{th mainrsult}, we need to prove the
inequality $v \le w^\sharp$.

   \begin{prop} \label{pr w le ksharp} $v\le w^{\sharp}$ on $[0,T]\x \Rpp\x [0,1]$.
      \end{prop}

\proof We use a backward induction argument.
We assume that \Hip~holds and that $v= w^{\sharp}$ and  on $[\tip,T]\x \Rpp \x [0,1]$ for some $i\le n-1$. Since it is true for $i=n$ by construction, the proof will be completed if one can show that   this implies that  \Hi~holds and that $v=w^{\sharp}$ on $[\ti,T]\x \Rpp \x [0,1]$.

Let us fix $(t,x,p)\in [\ti,\tip)\x \Rpp \x [0,1]$. Then, our induction hypothesis   implies that
$$
\E^{\Q_{t,x}}\left[\co{v\veee \ell}\left(\tip,X^{t,x}_{\tip},P^{t,p,\alpha}_{\tip}\right)\right]= \E^{\Q_{t,x}}\left[\co{w^{\sharp}\veee \ell}\left(\tip,X^{t,x}_{\tip},P^{t,p,\alpha}_{\tip}\right)\right]\,,
$$
for all $\alpha\in \Ac_{t,p}$. It then follows from Theorem \ref{kRepnew} and Proposition \ref{wProbRep1} that
$v(t,x,p)\le w^{\sharp}(t,x,p)$. But, the reverse inequality is proved in Proposition \ref{pr w ge ksharp}. This shows that $v=w^{\sharp}$ on $[\ti,T]\x \Rpp \x [0,1]$.
Then (i) of \Hi~is a consequence  of Proposition \ref{thm explicittildek}  and Proposition \ref{pr basic k and tilde k}. Proposition  \ref{wProbRep1} implies (ii) of \Hi.
Regarding the validity of (iii) of \Hi, it is proved in Lemma \ref{lem: induction Hip q-w} below.
 \eproof

\begin{lem}\label{lem: induction Hip q-w}  The hypothesis \Hip~implies \bru{{\rm\ref{item : Hip q-w}}} of \Hi.
 \end{lem}
 \proof
 It   follows from
\eqref{eq ktilde0} that
\begin{align} \label{eq expression tilde delta}
  q-(w^{\sharp}\veee \ell)^{\sharp}(t_{i},x,q)& =
\big(q-[q-\ell(t_{i},x)]^+\big){ \1_{A_1}}(t_i,x)
  \nonumber
\\
& +
(q-w(t_{i},x,q))\1_{A_2}(t_{i},x)
+\big(q-\kappa(t_{i},x,q)\big){\1_{A_3}}(t_i,x)\,,
\end{align}
in which
\begin{align*}
 q-\kappa(t_{i},x,q) =  (q - p_\ell(t_i,x)[q-q_\ell(t_i,x)]^{+}) \1_{\set{q \le
\bar{q}(t_i,x)}}
 + (q-w(t_{i},x,q))\1_{\set{q > \bar{q}(t_i,x)}} \,.
\end{align*}
By Lemma \ref{le more on ksharp}, $w^\sharp(t_i,x,1)=v(t_i,x,1)$ so that
$A_{2}\cup A_{3}=\set{v(\cdot,1)>\ell}$, recall \eqref{eq v(t,x,1) strict pos}.
In particular, we
observe that $\bar{q}<\infty$ on ${A_3}$.
The fact that the right-hand side in \eqref{eq expression tilde delta} converges
to $(v\veee \ell)(t_i,x,1)$ as $q\to \infty$ is then a consequence of  Lemma
\ref{le more on ksharp} and the definition of the $(A_i)_{i\le 3}$.


{It remains to show that each term in \eqref{eq expression tilde delta} is non-decreasing and continuous.  From   Lemma
\ref{le more on ksharp}, we know that $q\mapsto q-w(t_{i},x,q)$ is continuous and non-decreasing. The second term in the right-hand side of \eqref{eq expression tilde delta} is continuous and non-decreasing as well. As for the last term, we know that  $q\mapsto \kappa(\ti,x,q)$ is
continuous, so that it suffices to check the
monotony on each sub-interval   $(-\infty,\bar{q}(\ti,x)]$ and $[\bar
q(\ti,x),\infty)$ distinctly. On the second
interval,
we have that $q \mapsto q-\kappa(t_{i},x,q)$ is non-decreasing by Lemma \ref{le
more on ksharp}. This is also true on first interval since $p_\ell(t_i,x) \le
1$.}
\eproof

  \section{Appendix}

  We provide here the proofs of some technical results that were used in the proof of Theorem \ref{th mainrsult}.

   \vspace{2mm}

    {\bf Proof of Proposition \ref{prop : formulation sur temps d'adrets} } For $t=T$ the sets in  \eqref{eq: formulation Gamma avec temps arret} are $\R_{+}$ by definition of $\bT_{t}$ and ${\cal T}_{t}$. For $t<T$, {the definition of $\hat \tau_{\nu}$ implies
$
\S^{t,x,y,\nu}_{\hat \tau_{\nu}}{=} \bigcap_{s\in \bT_{t}} \S^{t,x,y,\nu}_{s}$, while
$ \bigcap_{s\in \bT_{t}} \S^{t,x,y,\nu}_{s}\subset \S^{t,x,y,\nu}_{\tau}
$,
 for any    $\tau\in {\cal T}_{t}$}. {Hence, for $t<T$, }
\begin{align*}
 { \p [\S^{t,x,y,\nu}_{\hat \tau_{\nu}} ]\ge p\;\Rightarrow\;\p\left[\bigcap_{s\in \bT_{t}} \S^{t,x,y,\nu}_{s}\right]\ge p
 \;\Rightarrow\;  {\p [ \S^{t,x,y,\nu}_{\tau} ]\ge
p,\;\forall\;\tau \in {\cal T}_{t}}\;\Rightarrow\;
 \p [\S^{t,x,y,\nu}_{\hat \tau_{\nu}} ]\ge p, }
 \end{align*}
 {where, in the last implication, we used the fact that $\hat \tau_{\nu} \in {\cal T}_{t}$.}  This proves   \eqref{eq: formulation Gamma avec temps arret} for $t<T$.
 \eproof

    \vspace{2mm}

{\bf Proof of Proposition \ref{wProbRep0}.} {\bf 1.} We first show that \eqref{eq probw0} holds.
Let $\bar v(t,x,p)$ denote the right-hand side of \eqref{eq probw0} and set
$$
J(t,x,p,\alpha):= \E^{\Q_{t,x}}\left[(v\veee \ell)\left(t_{i+1},X^{ t ,x}_{t_{i+1}},P^{t,p,\alpha}_{t_{i+1}}\right)\right]\,.
$$
Fix $y$ and $\alpha \in \Ac_{t,p}$ such that
$
y> J(t,x,p,\alpha).
$
Then, it follows from the martingale representation theorem that we can find $\nu  \in \Uc_{t,x,y}$ such that
$$
Y^{t,x,y ,\nu }_{\tip}> (v\veee \ell)\left(t_{i+1},X^{ t ,x}_{t_{i+1}},P^{t,p,\alpha }_{t_{i+1}}\right)\,.
$$
In particular, $Y^{t,x,y ,\nu }_{\tip}\ge \ell \left(t_{i+1},X^{ t ,x}_{t_{i+1}},P^{t,p,\alpha }_{t_{i+1}}\right)$. Since, we also have
 $Y^{t,x,y,\nu}_{\tip}$ $>$ $v (t_{i+1},X^{ t ,x}_{t_{i+1}},$ $P^{t,p,\alpha}_{t_{i+1}})$, it follows from the same arguments as in the proof of \cite[Lemma 2.2]{bouchard2010obstacle} that we can find a predictable process $(\tilde \nu,\tilde \alpha)$ which coincides with $(\nu,\alpha)$ on $[t,\tip]$, in the $\d t\x\d \P$-sense, and such that
 $$
 Y^{t,x,y,\tilde\nu}_{s}\ge \ell \left(s,X^{ t ,x}_{s},P^{t,p,\tilde\alpha}_{s}\right)\,, \mbox{ for all } s\in \bT_{\tip}\,.
 $$

 These processes are elements of $\hat \Uc_{t,x,y,p}$ whenever $\tilde \nu$ is
square integrable in the sense of \eqref{eq: cond carre int} and $\tilde
\alpha$ is such that $P^{t,p,\tilde \alpha}\in[0,1]$. The latter can be modified so that $P^{t,p,\tilde \alpha}$ is restricted to live in the interval
$[0,1]$ while $\tilde\nu$ can be modified so that \eqref{eq: borne Y} holds.
By the It\^{o} isometry, this induces the required square integrability property
of the financial strategy, recall \eqref{eq: cond mu sig}-\eqref{eq: cond
lambda}.
Combining the above with Proposition \ref{pr PbReduction0} shows that $\bar v(t,x,p)\ge v(t,x,p)$.

Conversely, let us fix $y>v(t,x,p)$. Then, it follows from the geometric dynamic programming principle of \cite[Theorem 2.1]{bouchard2010obstacle} that there exists  $(\nu,\alpha)\in \hat \Uc_{t,x,y,p}$ such that
$$
Y^{t,x,y,\nu}_{\tip} \ge (v\vee\ell)\left(\tip,X^{ t ,x}_{\tip},P^{t,p,\alpha}_{\tip}\right)\,.
$$
Since $Y^{t,x,y,\nu}$ is a super-martingale under $\Q_{t,x}$, this implies that $y\ge J(t,x,p,\alpha)$. The fact that $v(t,x,p)\ge \bar v(t,x,p)$ then follows from the arbitrariness of $\alpha$.

{\bf 2.} We now prove the Lipschitz continuity property. Note that it is true
for $t=T$, since $v(T,\cdot)=0$ by construction.
Let us assume that \eqref{eq: v loc lip} holds on $[\tip,T]$ for some $i<n$ and
show that it is then also true on $[\ti,T]$.
Let us fix $(t,p)\in [\ti,\tip)\times[0,1]$
and $x,x'\in \Rpp$.
We have that for all $\alpha\in \Ac_{t,p}$
\begin{align*}
&\left(Q^{t,x,1}_{\tip}\right)^{-1}(v\vee\ell)\left(\tip,X^{ t ,x}_{\tip},P^{t,p,\alpha}_{\tip}\right)
\\&\quad\quad\quad=\left(Q^{t,x',1}_{\tip}\right)^{-1}(v\vee\ell)\left(\tip,X^{ t ,x'}_{\tip},P^{t,p,\alpha}_{\tip}\right)
\\&\quad\quad\quad\quad+\left(Q^{t,x,1}_{\tip}\right)^{-1}\left[(v\vee\ell)\left(\tip,X^{ t ,x}_{\tip},P^{t,p,\alpha}_{\tip}\right)-(v\vee\ell)\left(\tip,X^{ t ,x'}_{\tip},P^{t,p,\alpha}_{\tip}\right)\right]
\\&\quad\quad\quad\quad+\left[\left(Q^{t,x,1}_{\tip}\right)^{-1}-\left(Q^{t,x',1}_{\tip}\right)^{-1}\right](v\vee\ell)\left(\tip,X^{ t ,x'}_{\tip},P^{t,p,\alpha}_{\tip}\right)\,.
 \end{align*}
Using first \eqref{eq: ass ell lipschitz}, the linear growth of $v$ (see \eqref{eq: borne v}) together with the fact that
\eqref{eq: v loc lip}  holds for  $(v\vee\ell)(\tip,\cdot,p)$, and using then \eqref{eq probw0}, we deduce that there exists $C>0$ such that $|v(t,x,p)-v(t,x',p)|$ is bounded by
$$
C \; \E^{\Q_{t,x}}\left[|X^{ t ,x}_{\tip}-X^{ t ,x'}_{\tip}|(1+ |X^{ t ,x}_{\tip}|+|X^{ t ,x'}_{\tip}|)+|Q^{t,x,1}_{\tip}/Q^{t,x',1}_{\tip} -1|\;(1+|X^{ t ,x'}_{\tip}|)\right]\,.
$$
In view of \eqref{eq: cond mu sig}-\eqref{eq: cond lambda}, this is controlled by $|x-x'|(1+|x|+|x'|)$ up to a multiplicative constant.
\eproof
     \vspace{2mm}

 {\bf Proof of Proposition \ref{prop: solution visco et unicite}.}  The growth property on $[0,T)\x \Rpp \x (0,\infty)$ follows from  Proposition \ref{pr basic k and tilde k} (which will be proved just below), Theorem \ref{th mainrsult}, \eqref{eq: extension v a R} and  \eqref{eq: borne v},
 $$
 0\le w(t,w,q)=\sup_{p\in \R} (pq-v(t,x,p))=\sup_{p\in [0,1]} (pq-v(t,x,p))\le q\,.
 $$
Note that Theorem \ref{th mainrsult} implies that $(w^{\sharp}\veee \ell)^{\sharp}(T,\cdot)=\ell^{\sharp}$. The fact that the lower- (resp. upper-) semicontinuous envelope of $w$ is a viscosity super- (resp. sub-) solution of ($\cal S$) is standard and we omit the proof. Continuity will then follow from the comparison principle. {The comparison can be  proved by backward induction. It is well-known that \eqref{eq ViscSol} admits a comparison principle in the class of functions with polynomial growth, see e.g.~\cite{CrIsLi92}. Hence, the comparison holds on $[t_{n-1},T)$.  Assume that it holds on $[\tip,T),\,i<n$ and that $(u_{j}^{\sharp}{\1_{[0,T)}}\veee \ell)^{\sharp}(\tip,\cdot)$ has polynomial growth, for $j=1,2$, then  it holds on $[\ti,T)$ too since $u_{1}(t_{i+1},\cdot)\ge u_{2}(t_{i+1},\cdot)$ implies $(u_{1}^{\sharp}\veee \ell)^{\sharp}(\tip,\cdot)\ge (u_{2}^{\sharp}\veee \ell)^{\sharp}(\tip,\cdot)$. Hence, we just have to prove that  $(u^{\sharp}_1\veee \ell)^{\sharp}$
has polynomial growth. By \cite[Theorem 16.5]{rocka}, we have $(u^{\sharp}_{j}\veee \ell)^{\sharp}=\co{u_{j}^{\sharp\sharp}\veee \ell^{\sharp}}$. Since $0\le u_{j}^{\sharp\sharp}\veee \ell^{\sharp} \le u_{j}\veee \ell^{\sharp}$ and the later has polynomial growth, the required property holds.}
 \eproof

          \vspace{2mm}

{\bf Proof of Proposition \ref{pr basic k and tilde k}.}   We proceed by
backward induction on $\bT_0 \cup \set{0}$.
Our claims are straightforward from
\eqref{eq: def w} at time $T$. Indeed, direct computations show that
$w^{\sharp}(T,\cdot,p)=0{+\infty\1_{\{p>1\}}}$.
Hence, $(w^{\sharp}\veee \ell)^{\sharp}(T,x,q)=
\ell^{\sharp}(T,x,q)=[q-\ell(T,x)]^{+}+\infty\1_{\set{q<0}}$. The properties
(a) and (b) hold.

We now assume that  (a) and (b) are satisfied on $[\tip,T]$ for some $i\le
n-1$ and fix $(t,x)\in [\ti,\tip)\x\Rpp$. Then, the definition of $w$ in
\eqref{eq: def w}   implies that $w(t,x,\cdot)$  is  non-negative,
non-decreasing, convex and that $w(t,x,0)=0$ (it is in particular proper).
It takes the value $+\infty$ for $q<0$, by \eqref{eq: def w} and  {the fact that
$(w^\sharp \vee \ell)^\sharp(t_{i+1},\cdot,q) = +\infty$ for $q <0$. Hence (a) holds on $[\ti,T]$.}  These two last assertions imply that $
{w}^\sharp(\cdot,p)=\sup_{q\ge0}\left\{pq-w(\cdot,q)\right\}$ and
$w^\sharp(t,\cdot,p) = 0$ for $p\le 0$.  We know from \cite[Theorem 12.2]{rocka}
that it is closed, convex and continuous on the interior of its domain. Since $w^\sharp$
is non-decreasing, by definition, we get from its closeness that it is
continuous on its domain.
The fact that $w^\sharp(t,\cdot,\cdot)\ge w^\sharp(t,\cdot,0) = 0$  also
implies that $( w^{\sharp}\veee\ell)(t,x,\cdot)$  is non-negative; moreover,
$(w^\sharp\veee\ell)(t,\cdot,0)=0$. For $q<0$, we then compute
$(w^{\sharp}\veee\ell)^{\sharp}(t,\cdot,q)= \sup_{p \le 1}\set{pq -
(w^{\sharp}\veee\ell)(t,\cdot,p)}= +\infty$. For $q \ge 0$, we get
$(w^{\sharp}\veee\ell)^{\sharp}(t,\cdot,q)= \sup_{p \in [0, 1]}\set{pq -
(w^{\sharp}\veee\ell)(t,\cdot,p)} \ge 0$. Moreover,
$(w^{\sharp}\veee\ell)^{\sharp}(t,x,\cdot)$ non-decreasing on $[0,\infty)$. By
definition, $(w^{\sharp}\veee\ell)^{\sharp}(t,x,\cdot)$ is closed, convex and
continuous on the interior of its domain. Being non-decreasing and closed, it
is in fact continuous on its domain.
  \eproof







\bibliographystyle{plain}

\end{document}